\documentclass[letterpaper, 10 pt, conference]{ieeeconf}

\newif\ifonlineapp %declare a plain version
\onlineapptrue
\onlineappfalse

\IEEEoverridecommandlockouts     % This command is only needed if you want to use the \thanks command

\usepackage{graphicx} % Required for inserting images

\makeatletter
\let\NAT@parse\undefined
\makeatother

\title{\LARGE \bf Asynchronous and Stochastic Distributed Resource Allocation}
\author{Qiang Li, Michal Yemini, and Hoi-To Wai% <-this % stops a space
\thanks{QL and HTW are with Department of SEEM, CUHK, Hong Kong SAR of China
        {\tt\small liqiang, htwai@se.cuhk.edu.hk}. MY is with Faculty of Engineering, Bar-Ilan University, Israel {\tt\small michal.yemini@biu.ac.il}.}%
}

\date{\today}

\usepackage[utf8]{inputenc} % allow utf-8 input
\usepackage[T1]{fontenc}
\usepackage{amsfonts} 
\usepackage{amssymb}

\usepackage{bm, Shortcuts_OPT}
\usepackage{amsthm}
\usepackage{amsmath}
\usepackage{url}
\usepackage{graphicx}

\allowdisplaybreaks[4]

\theoremstyle{plain}
\newtheorem{theorem}{Theorem}[section]

\theoremstyle{definition}

\theoremstyle{remark}

\newtheorem{lemma}{$\textrm{\bf Lemma}$\!}

\newtheorem{assumption}{$\textrm{\bf A}$\!\!}

\usepackage[colorlinks]{hyperref}
\hypersetup{
  colorlinks   = true, %Colours links instead of ugly boxes
  urlcolor     = blue, %Colour for external hyperlinks
  linkcolor    = red, %Colour of internal links
  citecolor    = red   %Colour of citations
}
\usepackage{cleveref}

\newcommand{\htwai}[1]{{\color[rgb]{1,0,0}[to: #1]}}

\begin{document}

\maketitle
\begin{abstract}
% \MY{This is an initial draft - nowhere near ready. I want to think a bit more about the flow}

This work proposes and studies the distributed resource allocation problem in \textit{asynchronous} and \textit{stochastic} settings. 
We consider a distributed system with multiple workers and a coordinating server with heterogeneous computation and communication times. 
We explore an approximate stochastic primal-dual approach with the aim of 1) adhering to the resource budget constraints, 2) allowing for the asynchronicity between the workers and the server, and 3) relying on the locally available stochastic gradients. 
We analyze our \textit{Asyn}chronous stochastic \textit{P}rimal-\textit{D}ual ({\algoname}) algorithm and prove its convergence in the second moment to the saddle point solution of the approximate problem at the rate of ${\cal O}(1/t)$, where $t$ is the iteration number. Furthermore, we verify our algorithm numerically to validate the analytically derived convergence results, and demonstrate the advantages of utilizing our asynchronous algorithm rather than deploying a synchronous algorithm where the server must wait until it gets update from all workers.

\end{abstract}

\begin{keywords}
    Distributed resource allocation, Asynchronous decentralized primal-dual algorithm.
\end{keywords}

\section{Introduction}
% \textcolor{red}{(MY) Should this text be part of the intro? If it where an ML venue possibly, but I'm not sure about a control one. Let's add some nice general discussion before it :)}

We consider the distributed resource allocation problem for a system with $n$ workers {and one coordinating server}, stated as the constrained optimization problem:
\beq\label{Q:orignal}
\begin{aligned}
\min_{\prm_i\in \RR^d, i\in [n]} &\quad f(\prm) \eqdef \sum_{i=1}^{n} f_i(\prm_i) 
\\
\text{subject to } &\quad \textstyle g_j\left( \sum_{i=1}^n A_{ji} \prm_i\right) \leq 0,~~j\in[m],
\\
&\quad \prm_i \in {\cal C}_i,~~ i\in[n].
\end{aligned}
\eeq 
For each worker $i \in [n]$, $\prm_i \in \RR^d$ denotes the $i$th local decision vector and  ${\cal C}_i$ is the feasible decision set. For each $j \in [m]$, $g_j: \RR^d \to \RR$ describes the (possibly non-linear) resource allocation constraint on a weighted sum of local decision vectors such that $A_{ji}$ is the weight on the $i$th decision vector {and $j$th (resource) constraint}. In the special case with $A_{ji} = 1/n$ for any $i,j$, the resource allocation constraint will be imposed on the averaged local decision vectors.
For brevity, we denote by $\prm \eqdef [\prm_1^\top, \prm_2^\top, \cdots, \prm_n^\top]^\top \in \RR^{nd}$ the joint decision vector.
Lastly, $f_i: \RR^d \to \RR$ is the local cost function.
In general, we have 
\beq\label{eq:stoc_loss}
    f_i(\prm_i) \eqdef \EE_{ Z_i\sim {\cal D}_i } [ \ell_i(\prm_i; Z_i)]
\eeq
where $Z_i \sim {\cal D}_i$ models the randomness or uncertainties in estimating the cost of a given decision.

Problem \eqref{Q:orignal} is motivated by multiple fields such as power distribution systems \cite{Pedram2010-smartGrid, fan2020online, li2011optimal, li2024socially}, wireless sensor networks with energy constraints \cite{ajay2022algorithm}, and traffic network congestion control \cite{kelly1998rate, Dong2023-network}.
For example, 
% {\color{blue} talk about \cite{ajay2022algorithm}}.
% For instance, 
\cite{kelly1998rate} analyzed the case of maximizing the network utility function given by $f_r(\prm_r) = \prm_{r} \log(\prm_r)$, where $\prm_r$ denotes the flow on route $r$ and the resource capacity constraints are formulated as $g_j(\prm_r) = {\sum_{r=1}^{R}} A_{jr} \prm_r - C_j \leq 0$, where $A_{jr} = 1$ if the resource $j$ is in route $r$, and $C_j$ is the capacity of resource $j$. {The work}
\cite{bakirtzis1994genetic} studied the economic dispatch problem with $f_i(\prm_i)$ representing the local generation cost function, subject to the inflow-outflow balance constraint $\sum_{i=1}^{n}\prm_i = \sum_{i=1}^{n} \prm_{d,i}$, where $\prm_{d,i}$ is the power output of generator $i$. {Additionally,}
\cite{li2011optimal, li2024socially} considered the power system utility maximization problem with $f_i( \prm_i )$ specifying the utility that user $i$ obtains and the resource allocation constraint $g_j(\cdot)$ refers to the cost for the utility company to provide $\sum_{i=1}^n \prm_i$ units of power to the customers. 
As noted by \cite{li2024socially}, the cost function $f_i(\prm_i)$ in the above examples can be subject to random fluctuations depending on the application scenarios, motivating the use of a stochastic cost model in \eqref{eq:stoc_loss}. {We refer the reader to \cite{doostmohammadian2025survey} for a comprehensive survey on applications.}
% Additionally, the generation limits are enforced through the constraint ${\cal C}_i = [\prm_i^{\sf min}, \prm_i^{\sf max}]$.
% Moreover, certain distributed robust optimization problems can  be cast in the form of \eqref{Q:orignal}. 
% For example, in \cite{duchi2019variance}, the author aim to minimize the empirical loss function $f(\prm)= \sum_{i=1}^{n} \prm_i z_{i}$ under the constraint ${\bm 1}^\top \prm = 1$. The optimization is performed over the $\rho$-neighborhood of the empirical distribution, defined as $\prm \in \left\{ \prm \in \RR_+^n: \frac{1}{2}\norm{n \prm -{\bm 1}}^2 \leq \rho\right\}$.

\begin{figure}[!t]
    \centering
    \begin{subfigure}[b]{1.0\columnwidth}
        \centering
        \includegraphics[width=\columnwidth]{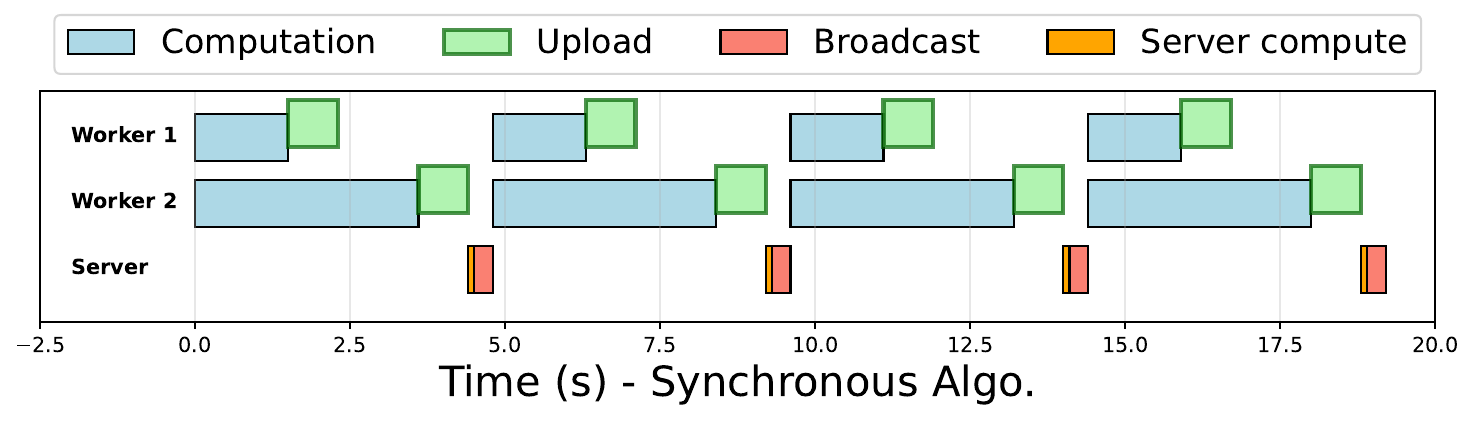}
        \caption{Synchronous Algorithm: Note that the system requires waiting for `Worker 1' before proceeding to the next step.}
        \label{fig:sync_algo}
    \end{subfigure}
    
    % \vspace{cm}
    
    \begin{subfigure}[b]{1.0\columnwidth}
        \centering
        \includegraphics[width=\columnwidth]{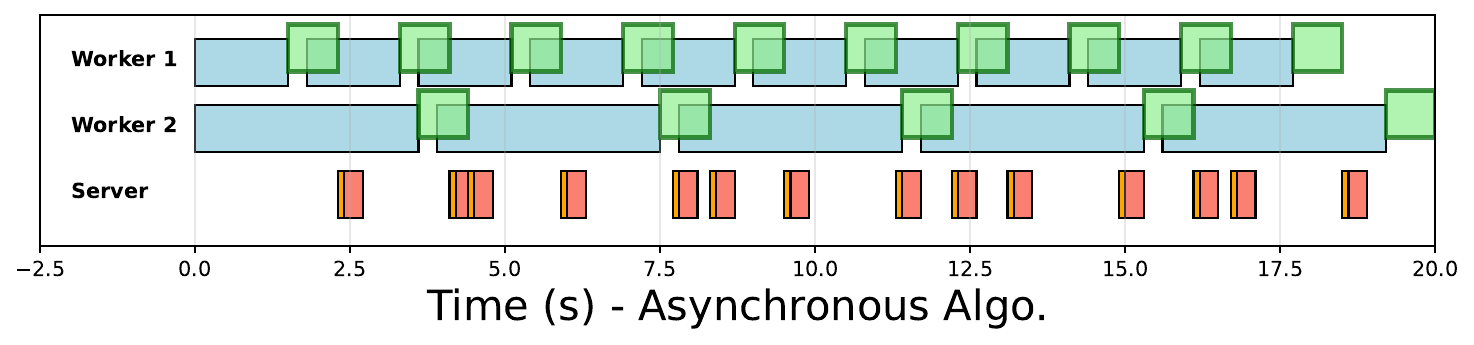}
        \caption{Asynchronous Algorithm: Note that the workers in the system are always active.}
        \label{fig:async_algo}
    \end{subfigure}
    \caption{Comparison of Synchronous and Asynchronous Algorithms for primal-dual algorithms tackling \eqref{Q:orignal}.}
    \label{fig:timeline_comparison}
    \vspace{-.8cm}
\end{figure}

Assuming a scenario where the workers are directly connected to a coordinating server, forming a star-like topology.
To develop distributed solutions for \eqref{Q:orignal}, a natural idea is to adopt the dual decomposition technique for handling the coupling constraints. This led to the regularized primal-dual algorithm proposed in \cite{koshal2011multiuser} which solves the deterministic of \eqref{Q:orignal}. Such algorithms are shown to have fast linear convergence while maintaining low computation complexity {for each worker}. Follow-up works have considered the robustification strategy against Byzantine attacks \cite{turan2020resilient}, online optimization with resource constraints \cite{fan2020online}, general forms of constrained optimization \cite{wu2022distributed}, etc.

A longstanding issue with distributed optimization is the \emph{straggler effects} caused by discrepancies in computation and communication speeds of heterogeneous workers. With a \emph{synchronous} algorithm, where the server waits until the updates of all workers arrive before executing its coordinating step, the straggler effect can lead to a significant slowdown to the system when $n \gg 1$, particularly if there is a worker that is significantly slower then the rest, as seen in Fig.~\ref{fig:timeline_comparison} (a).
An obvious remedy is to develop \emph{asynchronous} algorithms \cite{tsitsiklis1986} where the workers and server follow an independent schedule in updating their models; see Fig.~\ref{fig:timeline_comparison} (b).
Such architecture is relevant to power networks where asynchronous control is common \cite{wang2020asynchronous}. 
We note that for finite-sum optimization problems whose aim is to minimize $\sum_{i=1}^n f_i(\prm)$ in the absence of coupling constraints like in \eqref{Q:orignal}, asynchronous stochastic gradient (SG) algorithms have been extensively studied \cite{assran2020advances}. Examples include \cite{agarwal2011distributed} which studied the behavior of delayed SG methods, \cite{islamov2024asgrad, tyurin2023optimal} which developed tight bounds for variants of asynchronous SG methods, etc.
% koloskova2022sharper

For constrained optimization problems of the form \eqref{Q:orignal}, the study of asynchronous distributed algorithms is limited. Existing works require either synchronizing the dual variables periodically \cite{hale2017asynchronous}, or block-based parallel computations \cite{hendrickson2022totally}, and rely on deterministic gradients. 
This motivated the current paper to study: 
\begin{center}
    \emph{Can the problem \eqref{Q:orignal} be solved by an asynchronous \& stochastic distributed algorithm?}
    % \emph{Can the proposed {\algoname} algorithm solve problem \eqref{Q:orignal} in asynchronous setups? How does its solution quality compare to that of its synchronous counterpart?}
\end{center}  
We provide an affirmative answer to this question by proposing an asynchronous primal-dual ({\algoname}) algorithm which features two key aspects:  
1) A buffered central server that updates the dual variable upon receiving any worker's local model and broadcasts dual correction terms, which are stored in  local {worker's} buffers for future updates without interrupting workers.  
2) Buffered workers that operate independently, proceeding with local updates immediately after computation, ensuring that incoming server messages do not disrupt ongoing computations.  
Our main contributions are:

We show that the {\algoname} algorithm converges in expectation at a rate of ${\cal O}(1/t)$ with respect to the squared distance to the unique saddle point solution of a regularized version of \eqref{Q:orignal} under mild conditions, where $t$ is the number of iterations. The algorithm does not require a synchronized communication mechanism.  
    % \item To validate our theoretical findings, we conduct numerical simulations on a distributed Gaussian mean estimation problem with aggregated constraint budgets. We empirically demonstrate that our proposed method outperforms its synchronized counterpart.  
To our best knowledge, this is the first analysis to apply the asynchronous primal-dual algorithm to stochastic distributed resource allocation. 
% Our theoretical results provide strong evidence for the performance of {\algoname} compared to its synchronized counterpart. 

% asynchronous updates for straggler effects

% \noindent{\bf Related Works} \citep{turan2020resilient} investigates distributed algorithms based on primal-dual optimization in the presence of Byzantine attackers. Static and dynamic impersonation attack scenarios are analyzed, and corresponding robust optimization models are proposed to counter them. \citep{hale2017asynchronous} focuses on similar constrained problems and proposes an asynchronous multi-agent optimization framework in a fully decentralized setting, where each agent communicates only with its neighbors in the graph. They prove that the dual variable must remain synchronized across agents, while the primal update can be performed independently. 

% In \cite{even2024asynchronous}, the authors investigate unconstrained stochastic optimization and introduce an asynchronous SGD algorithm on graphs, where each worker performs local SGD and communicates with others to reach a consensus solution. \citep{wang2023linear} considers an incremental aggregated gradient method based on streaming data for a centralized worker-server architecture, where some gradients may be stale.
% \vspace{+.2cm}
% \noindent{{\bf Literature Review}. Shall we expand the asynchronized resource allocation literature?}

\vspace{+.2cm}
\noindent{\bf Notations.} For any $m \geq 1$, we denote $[m]\eqdef \{1, 2, \cdots, m\}$. $\EE[\cdot]$ denotes the full expectation over all randomness.
% Define the filtration $\mathcal{F}_k := \sigma \left(\prm_i^s, i \in[n], s=0, \ldots, k-1\right)$, where $\sigma(\cdot)$ denotes the sigma algebra generated by the random variables $Z$ in the operand. Observe that for any $k \geq$ 0 , $\prm^k$ is measurable with respect to (w.r.t.) $\mathcal{F}_k$. We shall use the shorthand notation $\mathbb{E}_k[\cdot]:=\mathbb{E}\left[\cdot \mid \mathcal{F}_k\right]$ for conditional expectation. 
$\mathds{1}(\cdot)$ is the indicator function. {Let $\supp({\cal D})$ be the support of distribution ${\cal D}$ and ${\sf conv}({\cal C})$ be the convex hull of set ${\cal C}$.}

\section{Problem Statement and Algorithm Design}\label{sec:algo}
This section reviews the synchronous primal-dual ({\synpd}) algorithm for the distributed resource allocation problem in \eqref{Q:orignal}. We then propose the \emph{asynchronous primal-dual} ({\algoname}) algorithm.
For simplicity, throughout this paper, we will consider a special case of \eqref{Q:orignal} such that $A_{ji} = 1/n$ for all $i,j$ such that the resource allocation constraint is imposed on the averaged decision vector.

\subsection{Preliminaries on Distributed Resource Allocation}
% Lagranian function -> add regularizer
    Observe that in \eqref{Q:orignal}, the objective function is decomposable while the only coupling between the local decision variables $\{ \prm_i \}_{i=1}^n$ lies in the constraints shaped by $g_j(\cdot)$, $j \in [m]$. To derive a distributed algorithm, we adopt the primal-dual framework. Introducing the dual variables $\blambda \in \RR_+^m$, under some regularity conditions, solving \eqref{Q:orignal} is equivalent to finding a saddle point to the following Lagrangian function: 
\begin{equation} \label{eq:lag_orig} \textstyle 
\sum_{i=1}^{n} f_i(\prm_i) +  \blambda^\top {\bm g}( \Bprm ),~~\Bprm := (1/n) \sum_{i=1}^n \prm_i,
\end{equation}
where ${\bm g}( \Bprm ) := [ g_1(\Bprm) ; \cdots ; g_m(\Bprm) ]$ collects the $m$  coupling resource allocation  constraints.

Similar to \cite{koshal2011multiuser}, we apply the dual smoothing technique through introducing the Tikhonov regularization to \eqref{eq:lag_orig}. Let $\upsilon>0$, we consider the regularized Lagrangian function:
% \MY{we use the term Lagrangian in the subsection title but don't explicitly use in the subsection itself.} 
\beq\label{def:lagrange}
    \textstyle {\cal L}(\prm, \blambda) \eqdef  \sum_{i=1}^{n} f_i(\prm_i) +       \blambda^\top {\bm g}(\Bprm) - \frac{ \upsilon }{2} \norm{\blambda}^2,
\eeq
% where $\Bprm \eqdef \frac{1}{n} \sum_{i=1}^{n} \prm_i$ is the aggregated model. $\blambda \in \RR^m_+$ is the dual variable, and $ (g_1, \cdots, g_m)^\top := {\bm g}: \RR^d \mapsto \RR^m$. 
Moreover, we consider the compact set $\bm{\Lambda} \subseteq \RR_+^m$ and the following min-max problem: 
\beq\label{Q:min-max}
\begin{aligned}
 \max_{\blambda\in \bm{\Lambda} } \min_{ \prm_i \in {\cal C}_i,  i\in [n]} &{\cal L}(\prm, \blambda).
\end{aligned}
\eeq
Notice that if each of $g_j(\cdot)$ is convex and each of $f_i$ is $\mu$-strongly-convex, then ${\cal L}(\cdot)$ is $(\mu,\upsilon)$-strongly convex-concave. Moreover, if $\blambda' \in \bm{\Lambda}$ (i.e., ${\rm diam}( \bm{\Lambda} )$ is sufficiently large), the unique saddle point solution to \eqref{Q:min-max} will be an {${\cal O}(\sqrt{\upsilon})$} approximation to the saddle point set of \eqref{eq:lag_orig}; see \cite[Lemma 3.3]{koshal2011multiuser}. 

To this end, a natural approach for solving \eqref{Q:min-max} is to adopt the primal-dual gradient descent-ascent method: for $k \geq 0$, 
\begin{subequations} \label{eq:pd_syn}
    \begin{align}
    \prm_i^{k+1} & = {\cal P}_{ {\cal C}_i } \left[ \prm_i^k - \gamma_{k+1} \grd_{ \prm_i } {\cal L}( \prm^k, \blambda^k ) \right],~i \in [n], \label{eq:pd_syn_a} \\
    \blambda^{k+1} & = {\cal P}_{ \bm{\Lambda} } \left[ \blambda^{k} + \gamma_{k+1} \grd_{ \blambda} {\cal L}( \prm^k, \blambda^k ) \right] \label{eq:pd_syn_b}
    \end{align}
\end{subequations}
Note that the primal-dual gradients can be computed as
\begin{align*}
     \grd_{\prm_i} {\cal L}(\prm^k, \blambda^k) &\eqdef \grd f_i(\prm_i^k) + \textstyle \frac{1}{n} \big( {\rm J} \, {\bm g}(\Bprm^k) \big) \blambda^k,~i \in [n],
    \\
    \textstyle \grd_{\blambda} {\cal L}(\prm^k, \blambda^k) &\eqdef {\bm g}(\Bprm^k) - \upsilon \blambda^k,
\end{align*}
where ${\rm J} \, {\bm g}(\Bprm) := [ \grd g_1(\Bprm)~\cdots~\grd g_m(\Bprm) ] \in \RR^{d \times m}$ is the Jacobian for the map ${\bm g}: \RR^d \to \RR^m$ at $\Bprm^k$. We denote the above algorithm \eqref{eq:pd_syn} as the synchronous primal-dual ({\synpd}) algorithm.

We notice that \eqref{eq:pd_syn} can be implemented via a distributed system. At any iteration $k \geq 0$, each worker receives the message containing $\frac{1}{n} \big( {\rm J} \, {\bm g}(\Bprm^k) \big) \blambda^{k}$ computed by the server and performs a local update following \eqref{eq:pd_syn_a}. Meanwhile, the server performs update for $\blambda^k$ \eqref{eq:pd_syn_b} based on the latest aggregated decision vector $\Bprm^k$. To prepare for the $k+1$th iteration, for the workers, the updated local decision vector $\prm_i^{k+1}$ has to be sent to the server and the server needs to broadcast $\frac{1}{n} \big( {\rm J} \, {\bm g}(\Bprm^{k+1}) \big) \blambda^{k+1}$ to the workers. We observe that this process is distributed where it does not require the $i$th worker to know the coupling constraints nor the local cost function of other workers, and vice versa. 

However, the above process requires \emph{synchronous} computations where the workers and server are not allowed to proceed to the $k+1$th iteration until the workers and server can respectively receive  $\frac{1}{n} \big( {\rm J} \, {\bm g}(\Bprm^{k+1}) \big) \blambda^{k+1}$, $\{ \prm_i^{k+1} \}_{i=1}^n$ from each other. 
If one (or more) of the workers is a \emph{straggler} such that it takes a significantly longer time to execute the gradient step \eqref{eq:pd_syn_a} and/or to send the updated model to the server, then the distributed process may be slowed down significantly. This issue is illustrated in Fig.~\ref{fig:timeline_comparison} (a) where the fast workers are left idle for a long time during the process, leading to inefficient use of computation resources.

\subsection{Asynchronous and Stochastic Primal-Dual Algorithm}
The straggler effect naturally leads us to consider an \emph{asynchronous update approach} where workers can continue their computations without waiting for others and the server's updates, e.g., see Fig.~\ref{fig:timeline_comparison} (b). To mimic \eqref{eq:pd_syn}, the server maintains a buffer that stores the last received decision vector from each worker and the update of $\blambda^k$ is performed using \emph{staled decision vectors}, subsequently the worker's local update will be affected by \emph{staled decision vectors} too. 

To formally describe and analyze the effects of asynchronous updates, we assume that the server's computation time is negligible compared to that of the workers' and the iteration index $k \geq 0$ counts the number of server's or workers' updates. We define 
\[
\begin{split}
    \tau_i^k \geq 0~ & \text{-- delay of $\prm_i$ stored at the server, $i \in [n]$,} \\
    \tau_s^k \geq 0~ & \text{-- broadcasting delay for $\big( {\rm J} \, {\bm g}(\Bprm) \big) \blambda$ from server.}
\end{split}
\]
That is, at the iteration $k$, the server has access to the aggregated and delayed decision vector $\Barb^k \eqdef \frac{1}{n} \sum_{i=1}^{n} \prm_i^{k-\tau_i^k}$, and the workers receive a delayed message given by ${\textstyle \frac{1}{n} \big( {\rm J} {\bm g}(\Barb^{k-\tau_s^k}) \big) \blambda^{k-\tau_s^k}}$.
Furthermore, 
% ${\cal A}^k$ is the set of active workers at iteration $k$, 
% ${\cal A}_i^k \eqdef \{i\in {\cal A}_k\}$ 
${\cal A}_i^k$ is the event where worker $i$ completes its computation at iteration $k$, and ${\cal E}_k$ is the event where the server receives at least one local decision vector at iteration $k$.

Define the sampled version of the Lagrangian function as:
\beq\label{def:lagrange_stoc}
    \textstyle L(\prm, \blambda; {\bm z}) \eqdef  \sum_{i=1}^{n} \ell_i(\prm_i; z_i) + \blambda^\top {\bm g}(\Bprm) - \frac{ \upsilon }{2} \norm{\blambda}^2,
\eeq
% \MY{Qiang, It seems that in Algorithms~\ref{algo:pd-server} \& \ref{algo:pd-worker}  both sides wait to receive input from the other - should we  explicitly mention of the initialization phase in this case? )}
We consider the asynchronous primal-dual ({\algoname}) algorithm for tackling \eqref{Q:min-max}:\vspace{-.3cm}
\begin{center}
\fbox{\begin{minipage}{.975\linewidth}\vspace{-.2cm}
\begin{align}
    &\prm_i^{k+1} = {\cal P}_{{\cal C}_i} \left[ \prm_i^k - \gamma_{k+1} \, \mathds{1}({\cal A}_i^k) \,  G_{i, \prm}^k(\prm_i^k , Z_i^{k+1} ) \right], \forall i \in [n], \notag
    \\
    &\blambda^{k+1} = {\cal P}_{\bm \Lambda}\left[ \blambda^{k} + \gamma_{k+1} \, \mathds{1}({\cal E}_{k}) \, G_{\blambda}^k(\blambda^k) \big) \right], \label{eq:asyn_pd}
\end{align}
\end{minipage}}
\end{center}
where $\mathds{1}( {\cal E} )$ is the indicator function that equals $1$ if the event ${\cal E}$ is true and is $0$ otherwise, $Z_i^{k+1} \sim {\cal D}_i$ is drawn independently. We also observe that 
\begin{align*}
    G_{i,\prm}^k ( \prm_i^k,  Z_i^{k+1} ) & \eqdef \grd \ell_i(\prm_i^k, Z_i^{k+1}) + \textstyle \frac{1}{n} \big( {\rm J} {\bm g}(\Barb^{k-\tau_s^k}) \big) \blambda^{k-\tau_s^k} ,
    \\
    G_{\blambda}^k (\blambda^k) & \eqdef {\bm g}(\Barb^{k}) - v \blambda^k,
\end{align*}
The above represent the delayed primal-dual stochastic gradients with respect to the primal and dual variables, and they approximate the primal-dual gradients for the sampled Lagrangian function $L( \prm^k, \blambda^k; {\bm Z}^{k+1} )$. The {\algoname} algorithm operates through two main components: the server update procedure and the worker update procedure, as detailed in Algorithms~\ref{algo:pd-server} \& \ref{algo:pd-worker}. Note that both routines are only triggered when messages are received from each side. 

\begin{algorithm}[t]
\caption{{\bf (Server)} {\algoname} algorithm}\label{algo:pd-server}
\begin{algorithmic}[1]
    \STATE {\bfseries Input}: Maximum iterations $K$, Step size rule $\gamma_{k} $.
    \STATE {\bfseries Initialize}: counter $k=0 $, central buffer ${\bm b}^{0}={\bm 0}$, dual variable $\blambda^0 = {\bm 0}$. 

    \WHILE{$k \leq K$}
            \STATE Wait until at least one local model copy $\prm_i^{k-\tau_{i}^k} $ sent by worker $i$ is received.
        \STATE Update the $i$-th slot of  buffer: ${\bm b}_{i}^{k} = \prm_i^{k-\tau_i^k} $. 
        \STATE Broadcast $\frac{1}{n}\big({\rm J} {\bm g} (\overline{\bm b}^k)\big) \blambda^{k} $ to all workers.
        \STATE Update the dual variable:\\
        $\blambda^{k+1} = {\cal P}_{\bf \Lambda} \left( \blambda^k + \gamma_{k+1} G_{\blambda}^{k} (\blambda^k) \right).$
        \STATE $k \leftarrow k+1$.
    \ENDWHILE
    \STATE %Boardcast signal EXIT and 
    Output dual variable $\blambda^{K+1}$.
\end{algorithmic}
\end{algorithm}

\begin{algorithm}[t]
\caption{{\bf (Worker $i$)} {\algoname} algorithm}\label{algo:pd-worker}
    \begin{algorithmic}[1]
        \STATE {\bfseries Input}: Maximum iterations $ K $, Step size rule $\gamma_{k} $.
        \STATE {\bfseries Initialize}: counter $k=0 $ and local model $ \prm_{i}^0 $

        \WHILE{$k\leq K$}
            \STATE Receive $\frac{1}{n}\big( {\rm J} {\bm g}(\overline{\bm b}^k) \big) \blambda^{k}$ from server and store it in the local buffer.
            \STATE Draw sample(s) $ Z_{i}^{k+1} \sim {\cal D}_{i} $.
            \STATE Update the local primal variable:
            
            $
            \prm_i^{k+1} = {\cal P}_{{\cal C}_i} \left( \prm_i^k - \gamma_{k+1}  G_{i, \prm}^{k}(\prm_i^k; Z_i^{k+1}) \right).
            $
            \STATE Send local model $\prm_i^{k+1}$ to the server.
            \STATE $k \leftarrow k+1$.
        \ENDWHILE
        \STATE Output local primal variable $\prm_i^{K+1}$.
    \end{algorithmic}
\end{algorithm}

As observed, the server updates the dual variable only when it receives a new local model copy, as indicated by $\mathds{1}({\cal E}_k)$. For each worker, $\mathds{1}({\cal A}_i^k)$ determines whether worker $i$ is active at iteration $k$. If active, the worker performs an update; otherwise, if the worker is occupied with computation, it essentially {\emph{retains the previous iterate}} $\prm_i^k$, which is crucial for convergence analysis.

Note that the {\algoname} algorithm involves delays in both directions of computation and communication delays between workers and server. 
{Different from the ASSP (Asynchronous Stochastic Saddle Point) algorithm studied by \cite{bedi2019asynchronous}, which primarily considers delays caused by communication channels, the {\algoname} algorithm operates in a fully asynchronous setting where delays arise not only from communication latency but also from heterogeneous and varying computation speeds across workers, making the asynchrony substantially more complex.} 
Compared to \cite{hale2017asynchronous}, the {\algoname} algorithm does not require a double loop design which resets the buffer variable periodically. Instead, the {\algoname} algorithm allows the workers and server to perform continuous updates. These challenges led us to develop a new convergence analysis in the next section.

\section{Convergence Analysis}

This section establishes the convergence of the {\algoname} algorithm to the saddle point $(\prm^\star, \blambda^\star)$ of \eqref{Q:min-max}. We first state the assumptions about the objective function, constraints, {stochasticity,} and the system delays. 

\begin{assumption}\label{assu:compact_set}
There exists constants $D, M>0$ such that for any $\prm_i, \prm_i^\prime \in {\cal C}_i$, $\blambda, \blambda^\prime \in {\bm \Lambda}$, ${\bm z} \in {\rm supp}( {\cal D} )^n$, it holds
\[
\begin{split}
& \textstyle \norm{\prm - \prm^\prime} \leq D ,~\norm{\blambda - \blambda^\prime} \leq D,~\max_{j\in [m]} \norm{\grd g_{j}(\Bprm)} \leq M, \\
& \textstyle \max_{i\in [n]} \norm{\grd_{\prm_i} {L}(\prm, \blambda; z)} \leq M,~\norm{\grd_{\blambda} {L}(\prm, \blambda ; z)} \leq M.
\end{split}
\]
    % \beqq
    %     \max\left\{ \max_{i\in [n]} \norm{\grd_{\prm_i} {L}_i(\prm_i, Z_i)}, \norm{\grd_{\blambda} {L}(\blambda)} \right\} \leq M
    %     \\
    %     \norm{\grd g_{j}(\Bprm)} \leq M,
    % \eeqq
    % where $Z_i\sim {\cal D}_i$, $\forall i \in [n]$, $\forall j \in [m]$.
\end{assumption}
The above assumption can be satisfied for bounded feasible sets ${\cal C}_i, \bm{\Lambda}$ together with the following conditions:
\begin{assumption}\label{assu:obj}
    For any $i\in[n]$, the loss function $f_i(\prm_i)$ is $L_i$-smooth such that for any $\prm^{\prime}, \prm \in {\cal C}_i$,
\beqq
    \left\|\nabla f_i(\prm)-\nabla f_i \left(\prm^{\prime}\right)\right\| \leq L_i \left\|\prm-\prm^{\prime}\right\|.
\eeqq
Moreover, there exists $\mu_i > 0$ such that
\beq \notag
    f_i(\prm') \geq f_i( \prm) + \pscal{ \grd f_i (\prm) }{ \prm' - \prm } + ({\mu_i}/{2}) \norm{ \prm' - \prm }^2.
\eeq 
\end{assumption}
\begin{assumption}\label{assu:cons}
For each $j \in [m]$, $g_j(\prm)$ is convex, and  {$L_j$-smooth such that} %its gradient is Lipschitz continuous, 
%i.e., 
for any $\Bprm^{\prime}, \Bprm \in {\rm conv}( {\cal C}_1, \ldots, {\cal C}_n )$,
    \beq
        \normtxt{\grd g_{j}(\Bprm) - \grd g_{j}(\Bprm^\prime)} \leq L_{j} \normtxt{\Bprm - \Bprm^\prime}.
    \eeq
\end{assumption}

Define the primal-dual gradient map for ${\cal L} ( \prm, \blambda )$ as:
\[ \Phi(\prm, \blambda) \eqdef 
\begin{pmatrix}
\grd_{\prm} {\cal L}(\prm, \blambda ) \\ 
- {\grd}_{\blambda} {\cal L}(\blambda, \blambda)
\end{pmatrix} \in \RR^{nd+m}.
\]
As a consequence of the above assumptions, the following lemma shows that $\Phi(\prm, \blambda)$ is a strongly monotone and Lipschitz continuous map.
\begin{lemma}\cite[Lemma 3.4]{koshal2011multiuser}\label{lem:grd_map}
Under A\ref{assu:compact_set}, \ref{assu:obj}, \ref{assu:cons}. 
For any ${\bm w} = (\prm, \blambda)$, ${\bm w}' = (\prm', \blambda')$, it holds that 
\beq 
\begin{split}
& \pscal{ \Phi( {\bm w} ) - \Phi( {\bm w}' ) }{ {\bm w} - {\bm w}' } \geq \mu \| {\bm w} - {\bm w}' \|^2 , \\
& \| \Phi( {\bm w} ) - \Phi( {\bm w}' ) \| \leq L \| {\bm w} - {\bm w}' \| 
\end{split}
\eeq 
where we have defined $\mu = \min_{i\in [n]} \{v, \mu_i\}$ and
% The regularized mapping $\Phi : \prod_{i=1}^{n} {\cal C}_i \times {\bm \Lambda} \mapsto \RR^{nd+m}$ is strongly monotone with constant $\mu = \min_{i\in [n]} \{v, \mu_i\}$ and Lipschitz with constant $L$ given by 
\[
\textstyle L = \sqrt{\left(L_{\sf max} + M + D L_{\sf g}\right)^2 + (M + \upsilon)^2}, L_{\sf g} = \sqrt{\sum_{j=1}^{m} L_j^2}.
\]
\end{lemma}

The following assumption is imposed on the variance of stochastic gradient:
\begin{assumption}\label{assu:var}
    For any $i\in [n]$ and fixed $\prm \in \mathbb{R}^d$, there exists $\sigma_i \geq 0$ such that
    \[
        \textstyle \mathbb{E}_{Z_i \sim \mathcal{D}_i}\left[\left\|\nabla \ell_i\left(\prm ; Z_i\right) \! - \! \nabla f_i(\prm)\right\|^2\right] \!\leq\! \sigma_i^2\left(1 \! +\! \left\|\prm-\prm^\star\right\|^2\right) .
    \]
\end{assumption}
Denote $\sigma_{\sf max}^2 \eqdef \max_{i\in [N]} \sigma_i^2$. Lastly, we assume that the delays introduced by asynchronous operations are bounded: 
\begin{assumption}\label{assu:delay}
There exists a positive integer $\overline{\tau} \in \mathbb{N}_{+}$ such that 
    \beq
    \textstyle \max_{i\in [n]} \left\{ \tau_i^k, \tau_s^k\right\} \leq \overline{\tau}.
    \eeq
\end{assumption}
\begin{assumption}\label{assu:exec}
There exists positive integers $p, B$ such that $1 \leq p\leq B$ and for any $k\geq 0$, it holds 
    \beq
        \sum_{s=k+1}^{k+B}    \mathds{1}(i \in {\cal A}_{s})  = p, ~\forall i \in [n], \quad \sum_{s=k+1}^{k+B}  \mathds{1}({\cal E}_s) = p.
    \eeq
\end{assumption}
This assumption ensures that within any time window of length $B$, each worker will be activated for $p$ times and the server will receive $p$ local updates. Note that $\bar{\tau} \leq B$ by design. We observe that $p/B$ is corresponds to the average update frequency of the asynchronous system. 

The convergence result for {\algoname} is summarized as:
\begin{theorem}\label{thm:asyn-spd}
    Under A\ref{assu:compact_set}--A\ref{assu:exec}. Suppose that the non-increasing step size $\{\gamma_{k}\}_{k\geq 1}$ satisfies $\forall k\geq 1$, $\ell\in[m]$,
    \[
        \sup_{k\geq 1} \gamma_{k} \leq \frac{2}{p \mu} , \quad \frac{\gamma_{k-(\ell +1)B+2}}{\gamma_{k - \ell B + 2}} \leq 1 + \frac{p\mu }{2} \gamma_{k-\ell B +2}.
    \]
    % Let $\w^k \eqdef \sum_{i=1}^{n} \normtxt{\prm_i^{k} - \prm_i^\star}^2 + \normtxt{\blambda^k - \blambda^\star}^2$. For any $k\geq 0$,  it holds 
    % \begin{align*}
    %     \EE \left[\w^{k+1}\right] \leq \prod_{i=1}^{m} \left( 1- p\mu \gamma_{k-iB+2} \right) \w^0  + \frac{2\C}{p \mu} \gamma_{k- B+2}.
    % \end{align*}
    % where $\C \eqdef B(n{\cal C}_1 + {\cal C}_2)+({\cal C}_3 + {\cal C}_{4})$. ${\cal C}_1, {\cal C}_2, {\cal C}_3, {\cal C}_4$ are constants defined in Lemmas \ref{lem:one_step_des}--\eqref{eq:def_C12}, \ref{lem:sumA}--\eqref{eq:def_C34}. Additionally, $m = \floor{(k+1)/{B}}$ and $\w^0$ denotes the initial error. 
    Let $\w^k \eqdef \sum_{i=1}^{n} \normtxt{\prm_i^{k} - \prm_i^\star}^2 + \normtxt{\blambda^k - \blambda^\star}^2$ denote the error metric at iteration $k$. For any $k\geq 0$, we have
    \begin{align}\label{eq:thm}
        \EE [\w^{k+1}] \leq \prod_{i=1}^{m} \left( 1- p\mu \gamma_{k-iB+2} \right) \w^0  \!+\! \frac{2\C}{p \mu} \gamma_{k- B+2},
    \end{align}
where $\C \eqdef B(n{\cal C}_1 + {\cal C}_2)+({\cal C}_3 + {\cal C}_{4})$ combines the constants defined in \eqref{eq:def_C12} and \eqref{eq:def_C34}, and $m = \floor{(k+1)/{B}}$ denotes the number of complete time windows.
\end{theorem}
The convergence rate in \eqref{eq:thm} is comparable to that of a {synchronous} stochastic primal-dual descent-ascent algorithm with Tikhonov regularization, e.g., \cite{shanbhag2013stochastic}. In particular, by setting $\gamma_t = c_0/ (c_1+t)$ with appropriate parameters $c_0,c_1>0$, \eqref{eq:thm} implies that $\EE [\w^{k+1}] = {\cal O}( nB/ (pk) )$.
Notice that the convergence rate in practice seems to be independent of that of the delay parameter $B/p$. 

% \begin{Corollary}
%     Under the same assumptions as in Theorem \ref{thm:asyn-spd}. When $p = B = 1$, it holds
%     \beq
%         \textstyle \EE[\Delta^{k+1}] \leq \prod_{i=1}^{k+1} (1-\mu \gamma_{i+1}) \Delta^0 + {\cal O}(\gamma_{k+1}/\mu).
%     \eeq
% \end{Corollary}

\subsection{Proof Outline}
This subsection outlines the proof of Theorem \ref{thm:asyn-spd}. To simplify notations, we define
\begin{align*}
    & \grd_{\prm_i}{\cal L}^\star \eqdef \grd_{\prm_i} {\cal L}(\prm^\star, \blambda^\star), \quad
    \grd_{\blambda}{\cal L}^\star \eqdef \grd_{\blambda} {\cal L}(\prm^\star, \blambda^\star), \\
    & \bm{\mathcal{G}}_{i,\prm}^k ( \prm_i^k ) \eqdef \grd f_i(\prm_i^k) + \textstyle \frac{1}{n} \big( {\rm J} {\bm g}(\Barb^{k-\tau_s^k}) \big) \blambda^{k-\tau_s^k}, \\
    & \bm{\mathcal{G}}^k_{\blambda} (\blambda^k) \eqdef {\bm g}(\Barb^{k}) - v \blambda^k, 
\end{align*}
where $(\{\prm_i^\star\}_{i=1}^n, \blambda^\star)$ is the unique saddle point of problem \eqref{Q:min-max}. Let $\Tprm^{k}_i \eqdef \prm_i^{k} - \prm_i^\star$, $\Tlambda^{k} \eqdef \blambda^k - \blambda^\star.$
We begin by observing the following asynchronous descent lemma:
\begin{lemma}\label{lem:one_step_des} Under A\ref{assu:compact_set}, \ref{assu:obj}, \ref{assu:var}. For any iteration $k \geq 0$, it holds
\begin{align*}
    &\EE \left[ \normtxt{\Tprm^{k+1}_{i}}^2 \right] \leq \EE\left[\normtxt{\Tprm_{i}^{k}}^2\right] + {\cal C}_1 \gamma_{k+1}^2 
    \\
    &\qquad \quad - \gamma_{k+1}\EE \left[ \mathds{1}({\cal A}^k_i) \pscal{\Tprm_i^k}{ \bm{\mathcal{G}}_{i,\prm}^k ( \prm_i^k )- \grd_{\prm_i} {\cal L}^\star } \right]
    \\
    &\EE\left[\normtxt{\Tlambda^{k+1}}^2\right] \leq \EE\left[\normtxt{\Tlambda^{k}}^2\right] + {\cal C}_2 \gamma_{k+1}^2 
    \\
    & \qquad \quad + \gamma_{k+1}\EE \left[ \mathds{1}({\cal E}_{k}) \pscal{\Tlambda^k}{ \bm{\mathcal{G}}^k_{\blambda} (\blambda^k) - \grd_{\blambda} {\cal L}^\star} \right]
\end{align*}
where ${\cal C}_1$, ${\cal C}_{2}$ are constants defined as
\begin{align}\label{eq:def_C12}
{\cal C}_1 &\eqdef 3 (\sigmamax^2 + (\sigmamax^2 + L_{\sf max}^2)D^2 + 6M^2D^2/n^2), \\
    {\cal C}_2 &\eqdef 2 D^2\left( M^2 + 2v^2 \right) .\notag
\end{align}
\end{lemma}
Summing the primal inequalities over $i=1,\ldots,n$, and incorporating the dual inequality, we obtain
\begin{align}
     & \EE \left[ \w^{k+1} \right] \leq \EE \left[ \w^k \right]  + \left(n {\cal C}_{1} + {\cal C}_{2} \right) \gamma_{k+1}^2
    \\
    &\quad  -  \gamma_{k+1} \EE \left[ \sum_{i=1}^{n} \mathds{1}({\cal A}^k_i)  \pscal{\Tprm_i^k}{ \bm{\mathcal{G}}^k_{i, \prm}(\prm_i^k) - \grd_{\prm_i} {\cal L}^\star } \right] \notag
    \\
    &\quad  - \gamma_{k+1} \EE \left[ \mathds{1}({\cal E}_{k}) \pscal{\Tlambda^k}{ \bm{\mathcal{G}}^k_{\blambda}(\blambda^k) - \grd_{\blambda} {\cal L}^\star} \right]. \notag 
\end{align}
Recall that $\w^k \eqdef \sum_{i=1}^{n} \normtxt{\prm_i^{k} - \prm_i^\star}^2 + \normtxt{\blambda^k - \blambda^\star}^2$. By adding and subtracting the non-delayed primal-dual gradients of ${\cal L}( \prm^k, \blambda^k )$, we have
\begin{align}\label{eq:w_des}
    \hspace{-.45cm} \EE [\w^{k+1}] \! \leq \! \EE[\w^{k}] \! - \! \gamma_{k+1} \EE [A_1^k + A_{2}^k]  \! + \! \left(n {\cal C}_{1} + {\cal C}_{2} \right) \gamma_{k+1}^2, 
\end{align}
where the terms $A_{1}^k, A_{2}^k$ are defined as 
% \htwai{please use the notation $\grd_{\prm_i} {\cal L}( \prm^k, \blambda^k )$ instead of $\grd_{\prm_i} {\cal L}_i( \prm_i^k )$ (also in the appendix proof), the latter has not been defined and it hides the dependence on $\blambda^k$.}
\begin{align*}
    A_{1}^{k} &\eqdef \sum_{i=1}^{n} \indA{i}{k} \pscal{\Tprm_i^k}{ \grd_{\prm_i} {\cal L}( \prm^k, \blambda^k ) - \grd_{\prm_i} {\cal L}^\star } 
    \\
    &\quad + \indE{k}\pscal{\Tlambda^k}{-\grd_{\blambda} {\cal L}(\prm^k, \blambda^k) + \grd_{\blambda} {\cal L}^\star},
    \\
    A_{2}^{k} &\eqdef \sum_{i=1}^{n}  \indA{i}{k} \pscal{\Tprm_i^k}{ \bm{\mathcal{G}}^k_{i, \prm}(\prm_i^k) - \grd_{\prm_i} {\cal L}( \prm^k, \blambda^k ) } 
    \\
    &\quad + \indE{k}\pscal{\Tlambda^k}{-\bm{\mathcal{G}}^k_{\blambda}(\blambda^k) + \grd_{\blambda} {\cal L}(\prm^k, \blambda^k)}.
\end{align*}
We have decomposed the inner product term by introducing the non-delayed gradient map $\Phi(\prm^k, \blambda^k)$, the discrepancy between the non-delayed and delayed gradients is captured in the term $A_2^k$.
Iteratively applying inequality \eqref{eq:w_des} over a time window $[k-B+1, B]$ yields:
\begin{align}
    \EE \left[\w^{k+1}\right] &\leq \textstyle \EE \left[ \w^{k-B+1} \right] + (n {\cal C}_1 + {\cal C}_2) \sum_{t=k-B+1}^{k} \gamma_{t+1}^2 \notag  
    \\
    &\quad \textstyle - \EE\left[ \sum_{t=k-B+1}^{k} \gamma_{t+1} \left( A_{1}^{t-1} + A_{2}^{t-1}\right)\right]. \notag
\end{align}
As $\{\gamma_{k}\}_{k\geq 1}$ is non-increasing sequence, we have
\begin{align}\label{eq:des_w_B}
    \EE \left[\w^{k+1}\right] &\leq \EE \left[ \w^{k-B+1} \right] + (n {\cal C}_1 + {\cal C}_2) B \gamma_{k-B+2}^2
    \\
    &\quad\textstyle - \EE\left[ \sum_{t=k-B+1}^{k} \gamma_{t+1} \left( A_{1}^{t-1} + A_{2}^{t-1}\right) \right]. \notag 
\end{align}
The above derivation prompts us to bound the sum of $A_1^t$ and $A_2^t$. A key observation is

\begin{lemma}\label{lem:sumA}(Inner product bound) Under A\ref{assu:compact_set}, \ref{assu:delay}, \ref{assu:exec}. For any iteration $k\geq 0$, it holds that 
\begin{align*}
   \textstyle  - \sum_{t=k-B+1}^{k} \gamma_{t+1} A_{1}^{t-1} &\!\geq\! -p \mu \gamma_{t-B+2} \w^{k-B+1} \!+ {\cal C}_3 \gamma_{t-B+2}^2, 
    \\
    \textstyle - \sum_{t=k-B+1}^{k} \gamma_{t+1} A_{2}^{t-1} &\!\leq\! {\cal C}_4 \gamma_{t-B+2}^2 
\end{align*}
where ${\cal C}_3, {\cal C}_4$ are constants, defined as following
\begin{align}
    {\cal C}_3 &\eqdef  2p(n+1)M \left( b^2 D + BDL + BM \right), \notag \\
    {\cal C}_4 &\eqdef {2p {\bar{\tau}} DM (M+1)} 
    % = 2p BDM (M+1)
    .  \label{eq:def_C34}
\end{align}
\end{lemma}
To leverage the smoothness and monotonicity of the gradient map $\Phi(\prm, \blambda)$ (Lemma \ref{lem:grd_map}), we adjust the indices of the local models and gradients, decomposing the sums in $A_1^t$ and $A_2^t$ into three terms and deriving bounds for each. Lemma \ref{lem:sumA} then provides a lower bound for the inner product term in \eqref{eq:des_w_B}:
\begin{align*}
    & \textstyle  - \sum_{t=k-B+1}^{k} \gamma_{t+1} (A_1^{t-1}+ A_{2}^{t-1}) 
    \\
    &\geq - \mu p \gamma_{t-B+2} \Delta^{k-B+1} +  {\cal O} \left( \gamma_{t-B+2}^2 \right).
\end{align*} 

\begin{proof}[\bf Proof of Theorem \ref{thm:asyn-spd}]
    
Combining \eqref{eq:des_w_B} and Lemma \ref{lem:sumA} leads to 
\begin{align*}
    \EE \left[\w^{k+1}\right] &\leq \left( 1 - p\mu \gamma_{k-B+2}\right) \EE[\w^{k-B+1}] + \C \gamma_{k-B+2}^2
\end{align*}
where $\C \eqdef B(n{\cal C}_1 + {\cal C}_2)+({\cal C}_3 + {\cal C}_{4}) $. Unrolling above recursion gives us
\begin{align*}
    \EE \left[\w^{k+1}\right] & \textstyle \leq \prod_{i=1}^{q} \left( 1 - p\mu \gamma_{k-iB+2}\right) \w^0 
    \\
    &~~ \textstyle + \C \sum_{i=1}^{m} \gamma_{k-iB+2}^2 \prod_{j=1}^{i-1} \left( 1 - p\mu \gamma_{k-jB+2}\right)
\end{align*}
where $q = \floor{\frac{k+1}{B}}$.
% $r = (k+1) {\sf mod } B$. 
$\w^0$ denotes the initial error. For edge cases, we denote that $\gamma_{k-iB+2} = \gamma_{1}$ and $\w^{k-iB+1} = \w^0$ if $k-iB+2 \leq 0$. Applying a variant of Lemma \ref{lem:aux} and requiring that for any $\ell\in [q]$,
\begin{align*}
   \textstyle  \sup_{k\geq 1}\gamma_{k} \leq \frac{2}{p \mu}, 
    \frac{\gamma_{k-(\ell +1)B+2}}{\gamma_{k - \ell B + 2}} \leq 1 + \frac{p\mu }{2} \gamma_{k-\ell B +2}, 
\end{align*}
it holds that
\begin{align*}
    \textstyle \EE \left[\w^{k+1}\right] \leq \prod_{i=1}^{q} \left( 1- p\mu \gamma_{k-iB+2} \right) \w^0  + \frac{2\C}{p \mu} \gamma_{k- B+2}.
\end{align*}
This completes the proof of Theorem \ref{thm:asyn-spd}.
\end{proof}

The omitted proofs in this section can be found in the online appendix \url{https://www1.se.cuhk.edu.hk/~htwai/pdf/cdc25-pd.pdf}.

\section{Numerical Simulation}\label{sec:simu}

This section presents simulation examples to validate our theoretical findings. All experiments are conducted with Python on a server using an Intel Xeon 6318 CPU.

\begin{figure}[!t]
    \centering    \includegraphics[width=1.\columnwidth]{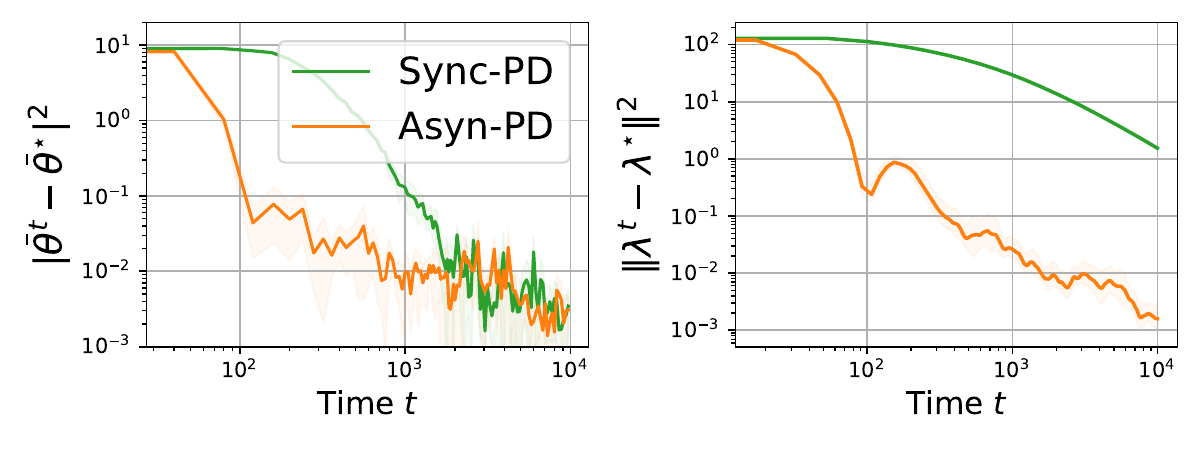}
    \caption{Primal and dual convergence comparison of Synchronous PD algorithm and {\algoname} algorithm. The shadow region represents the 90\% confidence interval. }
    % \htwai{Do we actually have the confidence interval?}
    \label{fig:simu}
    \vspace{-.4cm}
\end{figure}
\begin{figure}[!t]
    \centering    \includegraphics[width=1.\columnwidth]{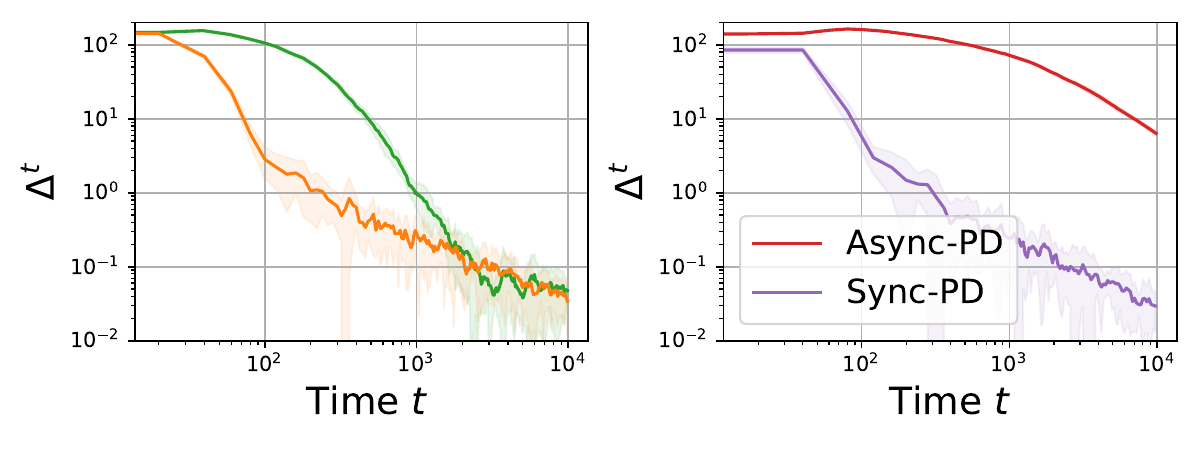}
    \caption{Comparison of the metric $\Delta^{t}$ between the Synchronous PD algorithm and the {\algoname} algorithm. \emph{(Left)} Simulation where each worker’s computation speed is set to $v = [4,4,3,2,1]$ per update. \emph{(Right)} An extreme case with $v = [10,4,3,2,1]$.}
    \label{fig:simu2}
    \vspace{-.4cm}
\end{figure}

% \vspace{+.2cm}

\noindent{\bf Resource Allocation Problem.} We consider problem \eqref{Q:orignal} with a worker-server architecture consisting of $n=5$ workers. The $i$th local stochastic objective function is given by $\ell_i(\prm_i; Z_i) = (\prm_i - Z_i)^2$, where the distribution for $Z_i$ is ${\cal D}_i \eqdef {\cal N}(\Bar{z}_i, \sigma^2)$, i.e., a normal distribution with mean $\Bar{z}_i = 10$ for $i=1,2,3$, $\Bar{z}_i = 12$ for $i=4, 5$ and standard deviation $\sigma=2$. Notice that this corresponds to a scenario with $f_i(\cdot) \neq f_j(\cdot)$, a.k.a.~a heterogeneous data setting.

% \htwai{Remember that we use a {\bf server-worker} architecture. (There's nothing like `learner')}
We set $m = 1$, and define the constraint function as $g(\mathbf{\Bprm}) \eqdef (1/n)\sum_{i=1}^n \prm_i - 5$. The feasible sets are $\mathcal{C}_i = [0, 7]$ for $i \in \{1, 2, 3\}$ and $\mathcal{C}_i = [0, 10]$ for $i \in \{4, 5\}$. The dual variable is constrained within ${\bm \Lambda} = [0, 10]$. The local worker models $\{\prm_i\}_{i=1}^{n}$ are randomly initialized and synchronized with the server’s buffer. The computational speed of each worker is randomly assigned as $v = [4, 4, 3, 2, 1]$ per update. The regularizer is set as $\upsilon=10^{-5}.$

% \htwai{$\upsilon=?$ -- I mean the regularizer! Also, $\upsilon \neq v$!}
For the given setting, the optimal solution is found to be as $\prm^\star = [4.26, 4.21, 4.25, 6.24, 6.16]$ and $\blambda^\star = 11.62$. For {\algoname}, the step size is set as $\gamma_{t} = {a_0}/{(a_1 + t)}$, with $a_0 = 10$ and $a_1 = 100$. 
% To validate the theoretical results, we implement {\algoname} in a sequential manner, allowing manual control of upload and broadcast delays to simulate a stable communication network. 
The upload delay is set as $\tau_i = 2$, $\forall i \in [n]$, and the broadcast delay as $\tau_s = 1$.  
As benchmark, we consider a \emph{synchronous PD} ({\bf Sync-PD}) algorithm where the server only performs an update for $\blambda^k$ after receiving the updated models from all workers, see Fig.~\ref{fig:timeline_comparison}. 

We run the stochastic algorithms for 10 repetitions and compare {\algoname} with {\bf Sync-PD}, reporting the average performances in Fig. \ref{fig:simu} in terms of the convergence of $\Bprm^k$, $\blambda^k$ towards $\prm^\star$, $\blambda^\star$, respectively. As observed, the {\algoname} algorithm converges to the optimal solution roughly at a rate of ${\cal O}(1/t)$, as predicted by Theorem \ref{thm:asyn-spd}. Our {\algoname} also converges faster than its synchronized counterpart to reach the same accuracy level due to the utilization of asynchronous updates.
In Fig.~\ref{fig:simu2}, we further consider a scenario with more heterogeneous computation speeds, denoted by $v' = [{10}, 4, 3, 2, 1]$ per update, i.e., the discrepancies in computation speed is more extreme. We compare the convergence in terms of the total error metric $\w^k$. The results demonstrate a more substantial performance degradation of {\bf Sync-PD} compared with {\algoname} when a single worker exhibits significantly slower computation speed.

\vspace{-.1cm}
\section{Conclusions}\vspace{-.1cm}
This study introduces an asynchronous primal-dual algorithm designed for stochastic distributed resource allocation. The proposed algorithm {\algoname} is shown to converge to an approximate primal-dual {saddle-point} solution at a rate of ${\cal O}({1/t})$. It exhibits better efficiency through fully utilizing the computational power of the distributed workers in the event of stragglers, compared to the {synchronous} {\synpd} algorithm. Numerical experiments corroborate our findings.

\vspace{-.2cm}

\bibliographystyle{ieeetr}
% Qiang: Let's use the palinnat format for convenience of citation and change it before we upload
\bibliography{ref.bib}

\begin{thebibliography}{10}

\bibitem{Pedram2010-smartGrid}
P.~Samadi, A.-H. Mohsenian-Rad, R.~Schober, V.~W.~S. Wong, and J.~Jatskevich,
  ``Optimal real-time pricing algorithm based on utility maximization for smart
  grid,'' in {\em 2010 First IEEE International Conference on Smart Grid
  Communications}, pp.~415--420, 2010.

\bibitem{fan2020online}
S.~Fan, G.~He, X.~Zhou, and M.~Cui, ``Online optimization for networked
  distributed energy resources with time-coupling constraints,'' {\em IEEE
  Transactions on Smart Grid}, vol.~12, no.~1, pp.~251--267, 2020.

\bibitem{li2011optimal}
N.~Li, L.~Chen, and S.~H. Low, ``Optimal demand response based on utility
  maximization in power networks,'' in {\em 2011 IEEE power and energy society
  general meeting}, pp.~1--8, IEEE, 2011.

\bibitem{li2024socially}
J.~Li, M.~Motoki, and B.~Zhang, ``Socially optimal energy usage via adaptive
  pricing,'' {\em Electric Power Systems Research}, vol.~235, p.~110640, 2024.

\bibitem{ajay2022algorithm}
P.~Ajay, B.~Nagaraj, and J.~Jaya, ``Algorithm for energy resource allocation
  and sensor-based clustering in m2m communication systems,'' {\em Wireless
  Communications and Mobile Computing}, vol.~2022, no.~1, p.~7815916, 2022.

\bibitem{kelly1998rate}
F.~P. Kelly, A.~K. Maulloo, and D.~K.~H. Tan, ``Rate control for communication
  networks: shadow prices, proportional fairness and stability,'' {\em Journal
  of the Operational Research society}, vol.~49, no.~3, pp.~237--252, 1998.

\bibitem{Dong2023-network}
S.~Dong, J.~Zhan, W.~Hu, A.~Mohajer, M.~Bavaghar, and A.~Mirzaei,
  ``Energy-efficient hierarchical resource allocation in uplink–downlink
  decoupled noma hetnets,'' {\em IEEE Transactions on Network and Service
  Management}, vol.~20, no.~3, pp.~3380--3395, 2023.

\bibitem{bakirtzis1994genetic}
A.~Bakirtzis, V.~Petridis, and S.~Kazarlis, ``Genetic algorithm solution to the
  economic dispatch problem,'' {\em IEE proceedings-generation, transmission
  and distribution}, vol.~141, no.~4, pp.~377--382, 1994.

\bibitem{doostmohammadian2025survey}
M.~Doostmohammadian, A.~Aghasi, M.~Pirani, E.~Nekouei, H.~Zarrabi, R.~Keypour,
  A.~I. Rikos, and K.~H. Johansson, ``Survey of distributed algorithms for
  resource allocation over multi-agent systems,'' {\em Annual Reviews in
  Control}, vol.~59, p.~100983, 2025.

\bibitem{koshal2011multiuser}
J.~Koshal, A.~Nedi{\'c}, and U.~V. Shanbhag, ``Multiuser optimization:
  Distributed algorithms and error analysis,'' {\em SIAM Journal on
  Optimization}, vol.~21, no.~3, pp.~1046--1081, 2011.

\bibitem{turan2020resilient}
B.~Turan, C.~A. Uribe, H.-T. Wai, and M.~Alizadeh, ``Resilient primal--dual
  optimization algorithms for distributed resource allocation,'' {\em IEEE
  Transactions on Control of Network Systems}, vol.~8, no.~1, pp.~282--294,
  2020.

\bibitem{wu2022distributed}
X.~Wu, H.~Wang, and J.~Lu, ``Distributed optimization with coupling
  constraints,'' {\em IEEE Transactions on Automatic Control}, vol.~68, no.~3,
  pp.~1847--1854, 2022.

\bibitem{tsitsiklis1986}
J.~Tsitsiklis, D.~Bertsekas, and M.~Athans, ``Distributed asynchronous
  deterministic and stochastic gradient optimization algorithms,'' {\em IEEE
  Transactions on Automatic Control}, vol.~31, pp.~803--812, Sept. 1986.

\bibitem{wang2020asynchronous}
Z.~Wang, F.~Liu, Y.~Su, P.~Yang, and B.~Qin, ``Asynchronous distributed voltage
  control in active distribution networks,'' {\em Automatica}, vol.~122,
  p.~109269, 2020.

\bibitem{assran2020advances}
M.~Assran, A.~Aytekin, H.~R. Feyzmahdavian, M.~Johansson, and M.~G. Rabbat,
  ``Advances in asynchronous parallel and distributed optimization,'' {\em
  Proceedings of the IEEE}, vol.~108, no.~11, pp.~2013--2031, 2020.

\bibitem{agarwal2011distributed}
A.~Agarwal and J.~C. Duchi, ``Distributed delayed stochastic optimization,''
  {\em Advances in neural information processing systems}, vol.~24, 2011.

\bibitem{islamov2024asgrad}
R.~Islamov, M.~Safaryan, and D.~Alistarh, ``Asgrad: A sharp unified analysis of
  asynchronous-sgd algorithms,'' in {\em International Conference on Artificial
  Intelligence and Statistics}, pp.~649--657, PMLR, 2024.

\bibitem{tyurin2023optimal}
A.~Tyurin and P.~Richt{\'a}rik, ``Optimal time complexities of parallel
  stochastic optimization methods under a fixed computation model,'' {\em
  Advances in Neural Information Processing Systems}, vol.~36,
  pp.~16515--16577, 2023.

\bibitem{hale2017asynchronous}
M.~T. Hale, A.~Nedi{\'c}, and M.~Egerstedt, ``Asynchronous multiagent
  primal-dual optimization,'' {\em IEEE Transactions on Automatic Control},
  vol.~62, no.~9, pp.~4421--4435, 2017.

\bibitem{hendrickson2022totally}
K.~R. Hendrickson and M.~T. Hale, ``Totally asynchronous primal-dual convex
  optimization in blocks,'' {\em IEEE Transactions on Control of Network
  Systems}, vol.~10, no.~1, pp.~454--466, 2022.

\bibitem{bedi2019asynchronous}
A.~S. Bedi, A.~Koppel, and K.~Rajawat, ``Asynchronous saddle point algorithm
  for stochastic optimization in heterogeneous networks,'' {\em IEEE
  Transactions on Signal Processing}, vol.~67, no.~7, pp.~1742--1757, 2019.

\bibitem{shanbhag2013stochastic}
U.~V. Shanbhag, ``Stochastic variational inequality problems: Applications,
  analysis, and algorithms,'' in {\em Theory driven by influential
  applications}, pp.~71--107, Informs, 2013.

\bibitem{li2022multi}
Q.~Li, C.-Y. Yau, and H.-T. Wai, ``Multi-agent performative prediction with
  greedy deployment and consensus seeking agents,'' {\em Advances in Neural
  Information Processing Systems}, vol.~35, pp.~38449--38460, 2022.

\end{thebibliography}

\newpage
\appendices

\section{An Auxiliary Lemma}
The following auxiliary lemma is quite standard, see \cite{li2022multi} Appendix E for detailed proof.

\begin{lemma}\label{lem:aux}
    Let $a>0$, $p\in \ZZ_+$ and $\left\{\gamma_{k}\right\}_{k \geq 1}$ be a non-increasing sequence such that $\gamma_{1}<2 / a$. If $\gamma_{k-1}^p / \gamma_{k}^p \leq 1+(a / 2) \gamma_{k}^p$ for any $k \geq 1$, then for any $k \geq 2$,
\[  
    \sum_{j=1}^{k} \gamma_{j}^{p+1} \prod_{\ell=1}^{j-1}\left(1-\gamma_{\ell} a\right) \leq \frac{2}{a} \gamma_{k}^p.
\]
\end{lemma}

\section{Proof of Lemma \ref{lem:one_step_des}}
\begin{proof}
    We observe the following chain according to the primal update rule,
    \begin{align*}
        &\norm{\prm_i^{k+1} - \prm_i^\star}^2 
        \\
        &\overset{(a)}{\leq} \norm{\prm_i^k - \gamma_{k+1} \indA{i}{k} {G}^k_{i, \prm}(\prm_i^k; Z_{i}^{k+1}) \! - \! \prm_i^\star \! + \! \gamma_{k+1}  \grd_{\prm_i} {\cal L}^\star)}^2
    \\
    &\leq \norm{\prm_i^k - \prm_i^\star}^2 + \gamma_{k+1}^2 \norm{    {G}^k_{i, \prm}(\prm_i^k; Z_{i}^{k+1}) -  \grd_{\prm_i} {\cal L}^\star  }^2
    \\
        &\quad - \gamma_{k+1} \indA{i}{k} \Pscal{\prm_i^k - \prm_i^\star}{   {G}^k_{i, \prm}(\prm_i^k; Z_{i}^{k+1}) -  \grd_{\prm_i} {\cal L}^\star }
    \end{align*}
where $(a)$ is due to the non-expansive property of projection operator. 
Taking expectation on both sides leads to
\begin{align}\label{eq:a}
    \EE[\normtxt{\Tprm_i^{k+1}}^2] &\leq \EE\left[\normtxt{\Tprm_i^k}^2\right]
    \\
    & + \gamma_{k+1}^2 \EE\left[ \normtxt{   G^k_{i, \prm}(\prm_i^k; Z_{i}^{k+1}) -  \grd_{\prm_i} {\cal L}^\star  }^2 \right] \notag 
    \\
    & \! - \! 2 \gamma_{k+1} \EE \! \left[ \! \indA{i}{k} \pscal{\Tprm_i^k}{  \bm{\mathcal{G}}^k_{i, \prm} (\prm_i^k) \! - \!  \grd_{\prm_i} {\cal L}^\star} \! \right].
    \notag 
\end{align}
Note $\Tprm^{k}_i \eqdef \prm_i^{k} - \prm_i^\star$. Next, let's deal with the second term,
\begin{align*}
    &\gamma_{k+1}^2 \EE \left[ \normtxt{  G^k_{i, \prm} (\prm_i^k;  Z_{i}^{k+1}) -  \grd_{\prm_i} {\cal L}^\star  }^2 \right] \\
    &\leq 3\gamma_{k+1}^2 \bigg( \EE \left[\norm{\grd \ell(\prm_i^k; Z_{i}^{k+1}) - \grd f_i(\prm_i^k)}^2 \right] 
    \\
        &\qquad + \normtxt{\grd f_i(\prm_i^k) - \grd f_i(\prm_i^k)}^2
        \\
        & \qquad  + \frac{1}{n^2} \normtxt{(\grd {\bm g}(\Barb^{k-\tau_s^k}))^\top \blambda^{k-\tau_s^k} - (\grd {\bm g}(\Bprm^\star))^\top 
        \blambda^\star}^2 \bigg)
    \\
    &\leq 3\gamma_{k+1}^2 \bigg( \sigmamax^2 (1+\norm{\prm_i^k - \prm_i^\star}^2) + L_i^2 \norm{\prm_i^k - \prm_i^\star}^2 
    \\
    &\qquad + \frac{1}{n^2} \normtxt{(\grd {\bm g}(\Barb^{k-\tau_s^k}))^\top \blambda^{k-\tau_s^k} - (\grd {\bm g}(\Bprm^\star))^\top 
    \blambda^\star}^2\bigg).
\end{align*}
For the last term, we observe the following bound,
\begin{align*}
    &\norm{(\grd {\bm g}(\Barb^{k-\tau_s^k}))^\top \blambda^{k-\tau_s^k} - (\grd {\bm g}(\Bprm^\star))^\top 
    \blambda^\star}^2 
    \\
    &\leq 2 \norm{(\grd {\bm g}(\Barb^{k-\tau_s^k}))^\top \blambda^{k-\tau_s^k} - (\grd {\bm g}(\Bprm^{k-\tau_s^k}))^\top 
    \blambda^\star}^2 
        \\
        &\qquad + 2\norm{(\grd {\bm g}(\Barb^{k-\tau_s^k}))^\top \blambda^\star - (\grd {\bm g}(\Bprm^\star))^\top 
        \blambda^\star}^2
    \\
    &\leq 2 \big\|\grd {\bm g}(\Barb^{k-\tau_s^k})\big\|^2 \normtxt{\blambda^{k-\tau_s^k} - \blambda^\star}^2 
    \\
        &\quad + 2\big\|\grd {\bm g}(\Barb^{k-\tau_s^k}) - \grd {\bm g}(\Bprm^\star)\big\|^2 \norm{\blambda^\star}^2 \leq 6 M^2 D^2.
\end{align*}
Substituting the above upper bound to \eqref{eq:a} gives us
\begin{align*}
    &\EE \norm{\Tprm_i^{k+1}}^2 \leq \EE \norm{\Tprm_i^k}^2 + {\cal C}_1 \gamma_{k+1}^2 
    \\
    &- 2 \gamma_{k+1} \indA{i}{k} \Pscal{\Tprm_i^k}{ \grd_{\prm_i} {\cal L}_i (\prm_i^k; \blambda^k) -  \grd_{\prm_i} {\cal L}^\star}
\end{align*}
where ${\cal C}_1 \eqdef 3 (\sigma_i^2 + (\sigma_i^2 + L_i^2)D^2 + 6M^2D^2/n^2)$. 

Next, we focus on the dual update rule, we have
\begin{align*}
    &\norm{\blambda^{k+1} - \blambda^\star}^2 
    \\
    &\leq \norm{\blambda^k - \blambda^\star + \gamma_{k+1} \left(  - G^k_{\blambda}(\blambda^k) + \grd_{\blambda}{\cal L}^\star \right)}^2
    \\
    &= \norm{\blambda^k - \blambda^\star}^2 + \gamma_{k+1}^2 \norm{ G^k_{\blambda}(\blambda^k) - \grd_{\blambda}{\cal L}^\star }^2 
    \\
    &\quad + \gamma_{k+1} \Pscal{\blambda^k - \blambda^\star}{- \grd_{\blambda}{\cal L}^\star + G^k_{\blambda}(\blambda^k) }.
\end{align*}
For the second term, we have
\begin{align*}
    &\gamma_{k+1}^2 \norm{ G^k_{\blambda}(\blambda^k) - \grd_{\blambda}{\cal L}^\star }^2 
    \\
    &= \gamma_{k+1}^2 \normtxt{ {\bm g}(\Barb^k) - v \blambda^k - {\bm g}(\Bprm^k) + v \blambda^\star}^2
    \\
    &\leq 2\gamma_{k+1}^2 \left( \norm{{\bm g}(\Barb^k) - {\bm g}(\Bprm^k)}^2 + v^2\norm{\blambda^k - \blambda^\star}^2\right)
    \\
    &\overset{(a)}{\leq}  2\gamma_{k+1}^2 \left( M^2 \normtxt{\Barb^k - \Bprm^k}^2+ 2v^2 D^2 \right)
    \\
    &\leq 2\gamma_{k+1}^2 \left( M^2 D^2 + 2v^2 D^2\right)
\end{align*}
where 
% in 
$(a)$
% , we use
% Lagrange's Mean Value Theorem and 
used A\ref{assu:cons}. Taking full expectation leads to 
\begin{align*}
    \EE \left[ \normtxt{\Tlambda^{k+1}}^2 \right] &\leq \EE\left[\normtxt{\Tlambda^k}^2\right] + 2\gamma_{k+1}^2 D^2\left( M^2 + 2v^2 \right) 
    \\
   & + \gamma_{k+1} \Pscal{\blambda^k - \blambda^\star}{ - \bm{\mathcal{G}}^k_{\blambda}(\blambda^k) + \grd_{\blambda}{\cal L}^\star}.
\end{align*}
This completes the proof of Lemma \ref{lem:one_step_des}.
\end{proof}

\section{Proof of Lemma \ref{lem:sumA}}
\begin{proof}
   \begin{figure*}[!t] 
    \hrulefill
    \begin{align}
        D_1 &\eqdef \sum_{i=1}^{n} \sum_{t \in S_i^{k-B+1: k}} \gamma_{t+1} \Pscal{\prm_i^{k-B+1} - \prm_i^\star}{\Tgrd_i {\cal L}^{k-B+1}} 
          + \sum_{t \in S_{n+1}^{k-B+1: k}} \gamma_{t+1} \Pscal{\blambda^{k-B+1} - \blambda^\star}{ \Tgrd_{\blambda} {\cal L}^{k-B+1}} 
        \notag\\
        D_2 &\eqdef \sum_{i=1}^{n} \sum_{t \in S_i^{k-B+1: k}} \gamma_{t+1} \Pscal{\prm_i^t - \prm_i^{k-B+1} }{\Tgrd_i {\cal L}^{t}} 
 + \sum_{t \in S_{n+1}^{k-B+1: k}} \gamma_{t+1} \Pscal{\blambda^{t} - \blambda^{k-B+1}}{ \Tgrd_{\blambda} {\cal L}^{t}} \label{eq:def_D123}
        \\
        D_3 &\eqdef \sum_{i=1}^{n} \! \sum_{t \in S_i^{k-B+1: k}} \! \! \!\!\!  \!\!\!  \gamma_{t+1} \Pscal{ \prm_i^{k-B+1} - \prm_i^\star}{\Tgrd_i {\cal L}^{t} - \Tgrd_{i} {\cal L}^{k-B+1}} 
           + \sum_{t \in S_{n+1}^{k-B+1: k}} \! \!\!\! \!\!\!  \! \gamma_{t+1} \Pscal{ \blambda^{k-B+1} - \blambda^\star}{\Tgrd_{\blambda} {\cal L}^{t} - \Tgrd_{\blambda}{\cal L}^{k-B+1}}
           \notag 
    \end{align}
    \hrulefill\vspace{-.5cm}
    \end{figure*}

We observe the following upper bound,
    \begin{align}\label{eq:sumA1}
        & \sum_{t=k-B+1}^{k} \gamma_{t+1} A_{1}^{t-1} 
        \\
        &= \sum_{t=k-B+1}^{k} \gamma_{t+1} \bigg( \sum_{i=1}^{n} \indA{i}{t} \pscal{\Tprm_i^t}{ \grd_{\prm_i} {\cal L}_i(\prm_i^t; \blambda^t) - \grd_{\prm_i} {\cal L}^\star } \notag 
        \\
        &\qquad + \indE{t}\pscal{\Tlambda^t}{-\grd_{\blambda} {\cal L}(\prm^t, \blambda^t) + \grd_{\blambda} {\cal L}^\star} \bigg) \notag 
        \\
        &\overset{(a)}{=}  \sum_{i=1}^{n} \sum_{t=k-B+1}^{k} \gamma_{t+1} \indA{i}{t} \cdot  \pscal{\Tprm_i^t}{ \grd_{\prm_i} {\cal L}_i (\prm_i^t, \blambda^t) - \grd_{\prm_i} {\cal L}^\star }  \notag 
        \\
            &\qquad + \sum_{t=k-B+1}^{k} \indE{t}\pscal{\Tlambda^t}{-\grd_{\blambda} {\cal L}(\prm^t, \blambda^t) + \grd_{\blambda} {\cal L}^\star} \notag
        \\
        & =: \left(\sum_{i=1}^{n} \sum_{t=k-B+1}^{k} \gamma_{t+1} \indA{i}{t} \cdot  \pscal{\Tprm_i^t}{ \Tgrd_{i} {\cal L}^t } \right)  \notag 
            \\ 
            &\qquad + \sum_{t=k-B+1}^{k} \indE{t}\pscal{\Tlambda^t}{ \Tgrd_{\blambda} {\cal L}^{t} } \notag 
    \end{align}
    where in equality $(a)$, we change the summation order. In the last equality, we use new notation of ${\Tgrd_i {\cal L}^t}$, ${\Tgrd_{\blambda} {\cal L}^t}$,
\begin{align*}
    {\Tgrd_i {\cal L}^t} &\eqdef \grd_{\prm_i} {\cal L}_i (\prm_i^t, \blambda^t) - \grd_{\prm_i} {\cal L}^\star,
    \\
    \Tgrd_\blambda {\cal L}^t &\eqdef -\grd_{\blambda} {\cal L}(\prm^t, \blambda^t) + \grd_{\blambda} {\cal L}^\star .
\end{align*}
Note that \eqref{eq:sumA1} has $n+1$ terms. We aim to analyze one of them, i.e., for any $i\in [n]$, it holds that 
\begin{align*}
        &\sum_{t=k-B+1}^{k} \gamma_{t+1} \indA{i}{t} \cdot  \pscal{\Tprm_i^t}{ \Tgrd_{i} {\cal L}^t } 
        \\
        &= \sum_{t=k-B+1}^{k} \gamma_{t+1} \indA{i}{t} \cdot  \pscal{\prm_i^t - \prm_i^\star}{ \Tgrd_{i} {\cal L}^t }
        \\
        &= \sum_{t=k-B+1}^{k} \gamma_{t+1} \indA{i}{t} 
            \bigg( \pscal{\prm_i^{k-B+1} - \prm_i^\star}{ \Tgrd_{i} {\cal L}^{k-B+1} } 
                \\
                &\qquad \qquad + \pscal{\prm_i^{t} - \prm_i^{k-B+1}}{ \Tgrd_{i} {\cal L}^t } 
                \\
                &\qquad\qquad  + \pscal{\prm_i^{k-B+1} - \prm_i^\star}{ \Tgrd_{i} {\cal L}^t - {\grd_i} {\cal L}^{k-B+1} } 
            \bigg)
        \\
        &= \sum_{t \in S_i^{k-B+1:k}} \gamma_{t+1} 
        \bigg( \pscal{\prm_i^{k-B+1} - \prm_i^\star}{ \Tgrd_{i} {\cal L}^{k-B+1} } 
            \\
            &\qquad \qquad + \pscal{\prm_i^{t} - \prm_i^{k-B+1}}{ \Tgrd_{i} {\cal L}^t } 
            \\
            &\qquad \qquad + \pscal{\prm_i^{k-B+1} - \prm_i^\star}{ \Tgrd_{i} {\cal L}^t - {\grd_i} {\cal L}^{k-B+1} } 
        \bigg)
    \end{align*}
    where we use add and subtraction trick to split one term to three terms. In the last equality, 
    \[
        S_i^{k-B+1: k} \eqdef \{t_i^1, t_i^2, \cdots, t_i^p\}, \quad i \in [n+1],
    \]
    where we let the server be the $(n+1)$th entity for simplicity, and $t_i^{j}$ is the $j$th time the indicator function $\indA{i}{j} = 1$. According to A\ref{assu:exec}, there are $p$ time points that worker $i$'s indicator function is 1. Correspondingly, we can rearrange R.H.S terms in \eqref{eq:sumA1} as following
    \begin{align*}
        \sum_{t=k-B+1}^{k} \gamma_{t+1} A_{1}^{t} =  D_1 + D_2 + D_3
    \end{align*}
    where $D_1, D_2, D_3$ are three terms defined in \eqref{eq:def_D123}.
 %    \begin{figure*}[!t] 
 %    \hrulefill
 %    \begin{align}
 %        D_1 &\eqdef \sum_{i=1}^{n} \sum_{t \in S_i^{k-B+1: k}} \gamma_{t+1} \Pscal{\prm_i^{k-B+1} - \prm_i^\star}{\Tgrd_i {\cal L}^{k-B+1}} 
 %          + \sum_{t \in S_{n+1}^{k-B+1: k}} \gamma_{t+1} \Pscal{\blambda^{k-B+1} - \blambda^\star}{ \Tgrd_{\blambda} {\cal L}^{k-B+1}} 
 %        \notag\\
 %        D_2 &\eqdef \sum_{i=1}^{n} \sum_{t \in S_i^{k-B+1: k}} \gamma_{t+1} \Pscal{\prm_i^t - \prm_i^{k-B+1} }{\Tgrd_i {\cal L}^{t}} 
 % + \sum_{t \in S_{n+1}^{k-B+1: k}} \gamma_{t+1} \Pscal{\blambda^{t} - \blambda^{k-B+1}}{ \Tgrd_{\blambda} {\cal L}^{t}} \label{eq:def_D123}
 %        \\
 %        D_3 &\eqdef \sum_{i=1}^{n} \! \sum_{t \in S_i^{k-B+1: k}} \! \! \!\!\!  \!\!\!  \gamma_{t+1} \Pscal{ \prm_i^{k-B+1} - \prm_i^\star}{\Tgrd_i {\cal L}^{t} - \Tgrd_{i} {\cal L}^{k-B+1}} 
 %           + \sum_{t \in S_{n+1}^{k-B+1: k}} \! \!\!\! \!\!\!  \! \gamma_{t+1} \Pscal{ \blambda^{k-B+1} - \blambda^\star}{\Tgrd_{\blambda} {\cal L}^{t} - \Tgrd_{\blambda}{\cal L}^{k-B+1}}
 %           \notag 
 %    \end{align}
 %    \hrulefill\vspace{-.5cm}
 %    \end{figure*}
    Next, we aim to obtain the upper bound for each term respectively. Since $D_1$ contains $\gamma_{t+1}$ which time index is varying, thus we further deal with it, i.e., 
    \[
        D_1 = D_{11} + D_{12},
    \]
    where $D_{11}$ and $D_{12}$ are defined as following,
    \begin{align*}
         D_{11} &\eqdef p \gamma_{k-B+2} \Pscal{\Tlambda^{k-B+1} }{ \Tgrd_{\blambda} {\cal L}^{k-B+1}}
        \\
        &+ \sum_{i=1}^{n} \sum_{t \in S_i^{k-B+1: k}} \gamma_{k-B+2} \Pscal{\Tprm_i^{k-B+1} }{\Tgrd_i {\cal L}^{k-B+1}} 
    \end{align*}
    \begin{align*}
        &D_{12} \eqdef 
        \\
        &\sum_{i=1}^{n} \sum_{t \in S_i^{k-B+1: k}}  \!(\!\gamma_{t+1} - \gamma_{k-B+2})  \Pscal{\Tprm_i^{k-B+1}}{\Tgrd_i {\cal L}^{k-B+1}} 
        \\
        & + \!\!\! \sum_{t \in S_{n+1}^{k-B+1: k}} (\gamma_{t+1} - \gamma_{k-B+2})  \cdot \Pscal{\Tlambda^{k-B+1}}{ \Tgrd_{\blambda} {\cal L}^{k-B+1}}
    \end{align*}
    For $D_{11}$, we can apply the strong monotonicity property of the gradient map $\grd {\cal L}(\prm, \blambda)$,
    \begin{align*}
        D_{11} &\geq \mu p \gamma_{k-B+2} \bigg( \sum_{i=1}^{n} \norm{\prm_i^{k-B+1} - \prm_i^\star}^2 
        \\
        &\quad + \norm{\blambda^{k-B+1} - \blambda^\star}^2\bigg) 
        = \mu p \gamma_{k-B+2} \w^{k-B+1}.
    \end{align*}
    We note that for constant step size $D_{12} = 0$. In general case, we know that $D_{12} \leq |D_{12}|$, i.e., 
    \begin{align*}
        & D_{12} \leq |D_{12}|
        \\
        % &\leq \sum_{i=1}^{n} \sum_{t \in S_i^{k-B+1: k}} \!(\gamma_{t+1} \!-\! \gamma_{k-B+2} ) \left| \Pscal{\prm_i^{k-B+1} \!- \! \prm_i^\star}{\Tgrd_i {\cal L}^{k-B+1}} \right| 
        % \\
        % & + \sum_{t \in S_{n+1}^{k-B+1: k}} (\gamma_{t+1} - \gamma_{k-B+2}) \left| \Pscal{\blambda^{k-B+1} - \blambda^\star}{ \Tgrd_{\blambda} {\cal L}^{k-B+1}} \right|
        % \\
        &\leq b^2 \gamma_{k-B+2}^2  \sum_{i=1}^{n}\sum_{t \in S_i^{k-B+1: k}} \norm{\Tprm_i^{k-B+1} } \cdot \norm{\Tgrd_i {\cal L}^{k-B+1}} 
        \\
        & + b^2 \gamma_{k-B+2}^2  \sum_{t \in S_{n+1}^{k-B+1: k}} \norm{\Tlambda^{k-B+1} } \cdot \norm{ \Tgrd_{\blambda} {\cal L}^{k-B+1}}
        \\
        &\overset{(a)}{\leq} b^2 \gamma_{k-B+2}^2 \cdot 2(n+1)p DM,
    \end{align*}
    where we use the fact that if $b \geq \frac{B}{a_0 (a_1 + t)}$, then $\gamma_t - \gamma_B \leq b^2 \gamma_{B}^2$. In $(a)$, we apply A\ref{assu:compact_set}. Therefore, 
    \begin{align}\label{eq:bnd_D1}
        D_1 \geq \mu p \gamma_{k-B+2} \w^{k-B+1} -   2b^2 (n+1)p D M \gamma_{k-B+2}^2 
    \end{align}
    For $D_2$, we have
\begin{align}\label{eq:D2}
    D_2 &= \sum_{i=1}^{n} \sum_{t \in S_i^{k-B+1: k}} \gamma_{t+1} \Pscal{\prm_i^t - \prm_i^{k-B+1} }{\Tgrd_i {\cal L}^{t}} \notag 
    \\
    & + \sum_{t \in S_{n+1}^{k-B+1: k}} \gamma_{t+1} \Pscal{\blambda^{t} - \blambda^{k-B+1}}{ \Tgrd_{\blambda} {\cal L}^{t}} 
\end{align}
We aim to analyze one of them, 
\begin{align*}
    &\sum_{t\in S_{i}^{k-B+1:k}} \gamma_{t+1} \Pscal{\prm_i^t - \prm_i^{k-B+1}}{\Tgrd_i {\cal L}^t} 
    \\
    &= \sum_{t\in S_{i}^{k-B+1:k}} \gamma_{t+1} \sum_{j=k-B+1}^{t-1} \Pscal{\prm_i^{j+1} - \prm_i^{j}}{\Tgrd_i {\cal L}^t}
    \\
    &\overset{(a)}{\leq} \sum_{t\in S_{i}^{k-B+1:k}} \gamma_{t+1} \sum_{j=k-B+1}^{t-1} \norm{\prm_i^{j+1} - \prm_i^{j}}\cdot \norm{\Tgrd_i {\cal L}^t}
    \\
    &\overset{(b)}{\leq}  \! \! \! \! \! \sum_{t\in S_{i}^{k-B+1:k}} \! \! \gamma_{t+1}  \! \sum_{j=k-B+1}^{t-1}  \! \norm{\gamma_{j+1} G^k_{i, \prm}(\prm_i^j; Z_{i}^{k+1})} \norm{\Tgrd_i {\cal L}^t}
    \\
    &\overset{(c)}{\leq} 2 M^2 \sum_{t\in S_{i}^{k-B+1:k}} \gamma_{t+1} \sum_{j=k-B+1}^{t-1} \gamma_{j+1} 
    \\
    &\leq 2p B M^2 \gamma_{t-B+2}^2 
\end{align*}
where in $(a)$, we use Cauchy-Schwarz inequality, in $(b)$ we apply the primal update, and $(c)$ is due to A \ref{assu:compact_set} \& \ref{assu:exec}. Substituting above upper bound back to \eqref{eq:D2} leads to
\begin{align}\label{eq:bnd_D2}
    D_2 \leq 2(n+1)p B M^2 \gamma_{t-B+2}^2 
\end{align}
Next, let's bound $D_3$,
\begin{align*}
    & D_3 = \sum_{i=1}^{n} \sum_{t \in S_i^{k-B+1: k}}  \! \! \! \! \gamma_{t+1} \Pscal{ \Tprm_i^{k-B+1} }{\Tgrd_i {\cal L}^{t} - \Tgrd_{i} {\cal L}^{k-B+1}} 
    \\
    &~~ + \sum_{t \in S_{n+1}^{k-B+1: k}} \gamma_{t+1} \Pscal{ \Tlambda^{k-B+1} }{\Tgrd_{\blambda} {\cal L}^{t} - \Tgrd_{\blambda}{\cal L}^{k-B+1}},
\end{align*}
we notice that there are $n+1$ terms, for each term, we have
\begin{align*}
    &\sum_{t \in S_i^{k-B+1: k}} \gamma_{t+1} \Pscal{ \Tprm_i^{k-B+1} }{\Tgrd_i {\cal L}^{t} - \Tgrd_{i} {\cal L}^{k-B+1} } 
    \\
    &= \sum_{t \in S_i^{k-B+1: k}} \gamma_{t+1} \sum_{j=k-B+1}^{t-1} \Pscal{ \Tprm_i^{k-B+1} }{\Tgrd_i {\cal L}^{j+1} - \Tgrd_{i} {\cal L}^{j} } 
    \notag
    \\
    &\leq D \sum_{t \in S_i^{k-B+1: k}} \gamma_{t+1} \sum_{j=k-B+1}^{t-1} \norm{\Tgrd_i {\cal L}^{j+1} - \Tgrd_{i} {\cal L}^{j} }
    \\
    &\overset{(a)}{\leq} D \sum_{t \in S_i^{k-B+1: k}} \gamma_{t+1} \sum_{j=k-B+1}^{t-1} 2 M L \gamma_{j+1}
    \\
    &\leq 2p BM DL \gamma_{k-B+2}^2
\end{align*}
where inequality $(a)$ is due to the fact that
\begin{align*}
   &\norm{\Tgrd_i {\cal L}^{j+1} - \Tgrd_{i} {\cal L}^{j} }
   \\
   &= \norm{\grd_{\prm_i} {\cal L}(\prm_i^{j+1}, \blambda^{j+1}) - \grd_{\prm_i} {\cal L}(\prm_i^{j}, \blambda^{j})} 
   \\
   &\overset{(b)}{\leq} L \left( \normtxt{\prm_i^{j+1} - \prm_i^j} + \normtxt{\blambda^{j+1} - \blambda^{j}}\right)
   \\
   &\leq L \left(  \gamma_{j+1}\norm{ G^k_{i, \prm} (\prm_i^j, Z_{i}^{j+1})} + \gamma_{j+1}\norm{ G^k_{\blambda}(\blambda^k)}\right) 
   \\
   &\leq 2M L \gamma_{j+1}. 
\end{align*}
where we use the smoothness property of the gradient map in $(b)$. Therefore, we obtain
\begin{align}\label{eq:bnd_D3}
    D_3 \leq 2(n+1)p BMDL \gamma_{k-B+2}^2.
\end{align}
Combining \eqref{eq:bnd_D1}, \eqref{eq:bnd_D2} and \eqref{eq:bnd_D3} gives us
\begin{align*}
    &- \sum_{t=k-B+1}^{k} \gamma_{t+1} A_{1}^{t-1} \leq - \mu p \gamma_{k-B+2} \w^{k-B+1} 
    \\
    &\qquad \quad + 
    2p(n+1)M \left( b^2 D + BDL + BM \right)\gamma_{k-B+2}^2.
\end{align*}
Then, we aim to bound $\sum_{t=k-B+1}^{k} \gamma_t A_{2}^{t-1}$. Recall that
\begin{align}
    &\textstyle \sum_{t=k-B+1}^{k} \gamma_{t+1} A_{2}^{t-1} = \sum_{t=k-B+1}^{k} \gamma_{t+1} \cdot \notag 
    \\
    & \quad \bigg[ \sum_{i=1}^{n}  \indA{i}{t} \pscal{\Tprm_i^t}{ \bm{\mathcal{G}}_{i, \prm}^t (\prm_i^t) - \grd_{\prm_i} {\cal L}(\prm_i^t, \blambda^t) } \label{eq:sumA2_1line}
    \\
    &\quad +  \indE{t}\pscal{\Tlambda^t}{- \bm{\mathcal{G}}_{\blambda}^t (\blambda^t) + \grd_{\blambda} {\cal L}(\prm^t, \blambda^t)} \bigg]. \label{eq:sumA2_2line}
\end{align}
For \eqref{eq:sumA2_1line}, we observe that 
\begin{align*}
    &\pscal{\Tprm_i^t}{ \bm{\mathcal{G}}^t_{i,\prm}(\prm_i^t) - \grd_{\prm_i} {\cal L}(\prm_i^t, \blambda^t) } 
    \\
    &\leq \norm{\Tprm_i^t} \cdot \normtxt{ \bm{\mathcal{G}}_{i,\prm}^t (\prm_i^t) - \grd_{\prm_i} {\cal L}(\prm_i^t) }
    \\
    &= \normtxt{\Tprm_i^t} \cdot \frac{1}{n} \norm{ (\grd {\bm g}(\Barb^{t-\tau_s^t}))^\top \blambda^{t-\tau_s^t} -( \grd {\bm g}(\Bprm^t))^\top \blambda^t  }
    \\
    &\leq \norm{\Tprm_i^t} \cdot \frac{1}{n} \bigg( \norm{ (\grd {\bm g}(\Barb^{t-\tau_s^t}))^\top \blambda^{t-\tau_s^t} -( \grd {\bm g}(\Bprm^{t-\tau_s^t}))^\top \blambda^t  } 
        \\
        &\qquad + \norm{ (\grd {\bm g}(\Barb^{t-\tau_s^t}))^\top \blambda^{t} -( \grd {\bm g}(\Bprm^t))^\top \blambda^t  } \bigg)
    \\
    &\leq \frac{D}{n} \left( M\normtxt{\blambda^{t-\tau_s^t} - \blambda^t} + D\normtxt{\grd {\bm g}(\Barb^{t-\tau_s^t}) - \grd {\bm g}(\Bprm^{t})} \right)
\end{align*}
where we apply A\ref{assu:compact_set}. Then, we focus on the difference between current $\blambda^k$ and its delayed copy, using A\ref{assu:delay}, we get
\begin{align}
    \norm{\blambda^{t -\tau_s^t} - \blambda^t} &\leq \sum_{j=t-\tau_s^t}^{t-1} \norm{\blambda^{j+1} - \blambda^j} 
    \\
    &= \sum_{j=t-\tau_s^t}^{t-1} \gamma_{j+1} \norm{ G_{\blambda}^j(\blambda^j)} \notag 
    \\
    &\leq \overline{\tau} M \gamma_{t-\bar{\tau} +2 } \leq B M \gamma_{k-B+2}. \notag 
\end{align}
\begin{align}\label{eq:diff-prm-b}
    \normtxt{\grd {\bm g}(\Barb^{t-\tau_s^t}) - \grd {\bm g}(\Bprm^{t})} &\leq M \normtxt{\Barb^{t-\tau_s^t} - \Bprm^t} 
    \\
    &\overset{(a)}{\leq} \frac{M}{n} \sum_{i=1}^{n} \norm{\prm_i^{t-\tau_s^t-\tau_i^{t-\tau_s^t}} - \prm_i^t} \notag 
    \\
    &\overset{(b)}{\leq} {2\bar{\tau} M^2\gamma_{k-\bar{\tau}+2}} 
    % 2B M^2\gamma_{k-B+2} 
    \notag 
\end{align}
where in inequality $(a)$, we recall that ${\Barb^{t}} = \frac{1}{n} \sum_{i=1}^{n} \prm_i^{t-\tau_i^t}$. In $(b)$, we use the following fact,
\begin{align*}
    \norm{\prm_i^{t-\hat{\tau}} - \prm_i^t} &\leq \sum_{j=t-\htau}^{t-1} \norm{\prm_i^{j+1} - \prm_i^{j}} 
    \\
    &\leq \sum_{j=t-\htau}^{k-1} \gamma_{j+1} \norm{ G_{i,\prm}^j(\prm_i^j; Z_{i}^{j+1})} 
    \\
    % &\leq 2B M\gamma_{k-B+2}, \\
    & {\leq 2 \bar{\tau} M \gamma_{k - \bar{\tau} + 1} }
\end{align*}
where we used $\hat{\tau} \eqdef \tau_s^t + \tau_i^{t-\tau_s^t} \leq 2 \bar{\tau}$ for simplicity. Combining above results, we have
\begin{align*}
    &\pscal{\Tprm_i^t}{ \bm{\mathcal{G}}_{i, \prm}(\prm_i^t) - \grd_{\prm_i} {\cal L}(\prm_i^t) }  
    {\leq \frac{\bar{\tau} D}{n} (M^2 + 2DM) \gamma_{k- \bar{\tau} +2}}.
    % \leq \frac{D}{n} (B M^2 + 2BDM) \gamma_{k-B+2} 
\end{align*}
Meanwhile, for \eqref{eq:sumA2_2line}, we have
\begin{align*}
    &\pscal{\Tlambda^t}{-\bm{\mathcal{G}}(\blambda^t) + \grd_{\blambda} {\cal L}(\blambda^t)}
    \overset{(a)}{\leq} \norm{\Tlambda^t} \cdot \norm{{\bm g}(\Barb^t) - {\bm g}(\Bprm^t)} 
    \\
    &\qquad \qquad \qquad \qquad \overset{(b)}{\leq} DM \norm{\Barb^t - \Bprm^t} {\overset{(c)}{\leq} \bar{\tau} DM^2 \gamma_{k-\bar{\tau}+2}}. 
    % \overset{(c)}{\leq} BDM^2 \gamma_{k-B+2}
\end{align*}
In inequality $(a)$, we use Cauchy-Schwarz inequality and the definition of $\grd_{\blambda} {\cal L}(\cdot)$. In the last inequality $(c)$, we apply same trick as in \eqref{eq:diff-prm-b}. Taking all results back,
\begin{align*}
    &\sum_{t=k-B+1}^{k} \gamma_{t+1} A_{2}^{t-1} 
    \\
    &\leq \sum_{t=k-B+1}^{k} \gamma_{t+1} \gamma_{k-B+2} \cdot {\bar{\tau}}
    \\
    &\quad \textstyle \cdot \bigg( \sum_{i=1}^{n} \indA{i}{t}  \frac{D}{n} (M^2 + 2DM)  + \indE{t} DM^2 \bigg)
    \\
    &\leq 2p {\bar{\tau}} DM (M+1) \gamma_{k-B+2}^2
\end{align*}
where in the last inequality, we use A\ref{assu:exec}.
\end{proof}

% {\color{red} TODO:
% \begin{itemize}
%     \item Hightlight the role of $\bar{\tau}$.
% \end{itemize}
% }

\end{document}